\documentclass[11pt]{article}
\usepackage{amsfonts, amsmath, amssymb, amsgen, amsthm, amscd,
latexsym}
\usepackage{color}
\usepackage[all]{xy}

\def\Id{\mathop{\rm Id}\nolimits}

\def\Hom{\mathop{\rm Hom}\nolimits}

\def\Id{\mathop{\rm Id}\nolimits}

\def\Cb{{\mathbb C}}

\def\Bc{{\cal B}}

\def\Fc{{\cal F}}

\def\Hc{{\cal H}}

\def\Kc{{\cal K}}

\def\Cc{{\cal C}}

\def\a{\alpha}

\def\d{\delta}
\def\D{\Delta}

\def\om{\omega}

\def\s{\sigma}

\def\t{\theta}

\def\ve{\varepsilon}

\def\ot{\otimes}

\def\ra{\rightarrow}

\def\rt{\triangleright}

\def\al{>\hspace{-4pt}\vartriangleleft}

\def\0D{\Delta^{(0)}}
\def\1D{\Delta^{(1)}}
\def\Db{\blacktriangledown}

\def\td{\tilde}

\newcommand{\Fs}{\mathfrak{s}}

\newtheorem{theorem}{Theorem}[section]
\newtheorem{remark}[theorem]{Remark}
\newtheorem{proposition}[theorem]{Proposition}
\newtheorem{lemma}[theorem]{Lemma}

\newtheorem{example}[theorem]{Example}
\newtheorem{definition}[theorem]{Definition}

\def\ni{\noindent}

\def\build#1_#2^#3{\mathrel{
\mathop{\kern 0pt#1}\limits_{#2}^{#3}}}
\newcommand{\ps}[1]{~\hspace{-4pt}_{^{(#1)}}}
\newcommand{\ns}[1]{~\hspace{-4pt}_{_{{<#1>}}}}

\numberwithin{equation}{section}
\def\a{\alpha}

\def\d{\delta}

\def\om{\omega}
\def\s{\sigma}
\def\t{\theta}
\def\ve{\varepsilon}

\def\D{\Delta}

\def\ot{\otimes}
\def\part{\partial}

\def\ra{\rightarrow}

\def\text{\hbox}

\def\ot{\otimes}

\def\ra{\rightarrow}

\def\Hom{\mathop{\rm Hom}\nolimits}

\def\Id{\mathop{\rm Id}\nolimits}

\def\build#1_#2^#3{\mathrel{
\mathop{\kern 0pt#1}\limits_{#2}^{#3}}}

\numberwithin{equation}{section}
\newcommand{\comment}[1]{\relax}

\begin{document}
%%%%%%%%%%%%%%%%%%%%%%%%%%%%%%%%%%%%%%%%%%%%%%%
\title{\bf \sc Equivariant Hopf Galois extensions and Hopf cyclic cohomology }

\author{
\begin{tabular}{cc}
Mohammad Hassanzadeh \thanks{Department of Mathematics  and   Statistics,
     University of New Brunswick, Fredericton, NB, Canada}\quad and \quad  Bahram Rangipour$~^\ast$
      \end{tabular}}

\maketitle

%%%%%%%%%%%%%%%%%%%%%%%%%%%%%%%%%%%
\begin{abstract}
\ni We define the notion of equivariant $\times$-Hopf Galois extension and apply it as a functor between the categories of SAYD modules of the $\times$-Hopf algebras involving in the extension. This generalizes the result of Jara-Stefan and  B\"ohm-Stefan on associating a SAYD modules to any ordinary Hopf Galois extension.
\end{abstract}

%%%%%%%%%%%%%%%%%%%%%%%%%%%%%%%%%%%%%%%%%%%%%%%%%%%
\section*{Introduction}
Hopf cyclic cohomology, which is now a well  known cohomology theory in noncommutative geometry,  was invented by A. Connes and H. Moscovici in \cite{ConMos:HopfCyc} as a computational tool for calculating the local index formula of  spectral triples associated to hypoecliptic operators on manifolds. Since then, the theory has been evolving to cover many other cases including algebras and coalgebras endowed with (co)symmetry from  Hopf algebras, bialgebras,  and  Hopf algebras with several objects such as para-Hopf algebras and $\times$-Hopf algebras \cite{ConMos:Cyc-Symm,ConMos:DiffCyc,bs2,Cra,Gor:SecChar,HaKhRaSo1,
HaKhRaSo2,JaSt:CycHom,Kay,Kay:UniHCyc,KhaPou,KhaRan:ANewCyc,KhaRan:ANoteCycDual,
KhaRan:Para_Hopf,
KhaRan:InvInvCycHom,Rang:CycCor}.  The main application of Hopf cyclic cohomology is to produce special cyclic cocycles on (co)algebras endowed with a (co)symmetry by a Hopf algebra. Generally speaking, Hopf cyclic cohomology of a Hopf algebra is more computable than cyclic cohomology of (co)algebras upon which the Hopf algebra (co)acts.
This accessibility is the main motivation and application of Hopf cyclic cohomology.
The Hopf cyclic cohomology of Hopf (co)module (co)algebras with coefficients in   SAYD modules was first defined in  \cite{HaKhRaSo2,HaKhRaSo1}, where  it was shown that the modular pair in involution defined by Connes and Moscovici  in \cite{ConMos:Cyc-Symm},  is in fact, a  one dimensional SAYD module. Later on,  Jara and Stefan  in \cite{JaSt:CycHom}  associated a canonical SAYD module to any Hopf-Galois extension such that the relative cyclic complex of the extension is isomorphic with  the cyclic complex defined on  the Hopf algebra with coefficients in the associated SAYD module.  Soon after, Khalkhali and the second  author showed in \cite{KhaRan:ANoteCycDual} that the  cyclic module obtained by Jara-Stefan is of the form of cyclic module defined in \cite{KhaRan:ANewCyc} which is  dual (in the sense of Connes' cyclic category) of Connes-Moscovici  cocyclic module.

 The first step for defining  Hopf cyclic cohomology of Hopf algebras with several objects was taken in \cite{ConMos:DiffCyc}, where it was shown, among other things, that  the symmetry on the crossed product algebra of diffeomorphism  group and algebra of smooth  functions on the frame bundle of a not necessarily flat manifold is defined by a Hopf algebroid generalizing Connes-Moscovici Hopf algebra in the case the manifold in question is flat.
 The cohomology defined in \cite{ConMos:DiffCyc} was abstractly generalized for a class of Hopf algebroids called Para Hopf algebras in \cite{KhaRan:Para_Hopf}. The Hopf cyclic cohomology of a coring with symmetry from a Para Hopf algebra defined in \cite{Rang:CycCor}. Finally Hopf cyclic cohomology of $\times$-Hopf algebras with coefficients in the SAYD modules  was defined  by B\"ohm and Stefan in \cite{bs2}, where it was also shown that the result of Jara-Stefan is valid for $\times$-Hopf Galois extensions.

 Let us  recall very briefly here the result  of Jara-Stefan from \cite{JaSt:CycHom}. Let  $\nobreak{\Db:A\ra A\ot H}$ define a comodule algebra and $B$ be the subalgebra of coinvariants for this coaction. Then the extension $B\subseteq A$ is called Galois if the canonical map ${ can}:A\ot_B A\ra  A\ot H$ is bijective. One uses the map $can$ iteratively to transfer the cyclic structure on the relative cyclic complex of $B\subset A$ to get a cyclic module on the other side depending only on the Hopf algebra $H$ and the SAYD module $A_B:= A/[A,B]$.

In this paper we start from the fact that relative cyclic homology of an extension $B\subseteq A$ is in fact Hopf cyclic homology of the ring $A$ over $B$  with coefficients in $B$ as a SAYD module over the $\times_B$-Hopf algebra $B^e:=B\ot B^{op}$. It is also easy to observe that $A_B= B\ot_{B^e} A$. On the other hand, if $M$ is a SAYD module over $B^e$, in the sense of \cite{bs2}, then the cyclic homology of $B$-ring $A$ with coefficients in $M$ under symmetry of  $B^e$ is well-defined. So  in case the extension is $H$-Galois   iterative application of the map $can$ takes us to a  cyclic complex which is again the dual of the Connes-Moscovici cocyclic module associated to the Hopf algebra $H$  with coefficients in $M\ot_{B^e} A$.
The mentioned  transfer of cyclic structures works  due to the fact that the map  $can$ happens to be $B^e$-equivariant. This means that, if one assumes the existence of a $\times$ Hopf algebra $\Kc$ acting on $A$  such that the map $can$ is equivariant then to any SAYD module $M$ over $\Kc$ one associates a SAYD module $\widetilde{M}:=M\ot_\Kc A$ over $H$. One then shows that the Hopf cyclic homology of the ring $A$ under symmetry of $\Kc$ with coefficients in $M$ is isomorphic with the dual cyclic module of Hopf cocyclic module of $H$ with coefficients in $\widetilde{M}$. The next step is to upgrade everything to the level of a $\times$-Hopf Galois extension, which is stated as Theorem \ref{general}.

The plan of the paper is as follows. In Section \ref{dual} we recall the basics of co(cyclic) modules and duality in the cyclic category. In Section \ref{S-2} we review the concepts of left and right $\times$-Hopf algebras and basics of  Hopf cyclic (co)homologies together with  some examples. In Section \ref{S-3} we define the equivariant $\times$-Hopf Galois extensions and prove the main result of the paper which is summarized in   Theorem \ref{general} and finally we bring some non trivial examples of this result.

Throughout the  paper we assume all objects are $\Cb$-vector spaces although everything works for $k$-modules, where $k$ is a commutative ring. We use the Sweedler  summations: for comultiplication of coalgebras or corings i.e., $\D(c)=c^{(1)}\ot c^{(2)}$, for  coactions i.e, $\Db(m)=m\ns{0}\ot m\ns{1}$; for translation map i.e,  $can^{-1}(b)=b\ns{-}\ot b\ns{+}$; and finally for the ``antipode" of $\times$-Hopf algebras, i.e, $\nu^{-1}(1\ot b)=b^-\ot b^+$.

\bigskip

{\bf Acknowledgement: }\\
The authors would like to thank the referee for his(her) carefully reading the manuscript and also for her(his) valuable comments.

%%%%%%%%%%%%%%%%%%%%%%%%%%%%%%%%%%%%%%%%%%%%%%%%%%%%%%%%%%%%%%%%%%%%%%%%%%%%%%%%%%%%%%%%%%%

 %%%%%%%%%
\tableofcontents

%%%%%%%%%%%%%%%%%%%%%%%%%%%%%%%%%%%%%%%%%%%%%%%%%%%%%%%%%%%%%%%%%%%%%%%%%%%%%%%%%%%%%%%%%%%%%%%%%%%%%%%%%%
\section{Basics of (co)cyclic modules}\label{dual}
In this section we recall the definitions of cyclic and cocyclic modules from \cite{NCG} (see also \cite{Loday}). Recall that a cosimplicial module is given by datum $(C_{n},\delta_{i}, \sigma_{i})$ where $C_{n}$, $n\geq 0$ is a $\mathbb{C}$-module. The maps $\delta_{i}:C^{n}\longrightarrow C^{n+1}$ are called cofaces, and $\sigma_{i}:C^{n}\longrightarrow C^{n-1}$ called codegeneracies. These are $\mathbb{C}$-module maps satisfying the following cosimplicial relations:
\begin{eqnarray}
\begin{split}
& \delta_{j}  \delta_{i} = \delta_{i} \delta_{j-1},   \hspace{35 pt}  \text{ if} \quad\quad i <j,\\
& \sigma_{j}  \sigma_{i} = \sigma_{i} \sigma_{j+1},    \hspace{30 pt}  \text{    if} \quad\quad i \leq j,\\
&\sigma_{j} \delta_{i} =   \label{rel1}
 \begin{cases}
\delta_{i} \sigma_{j-1},   \quad
 &\text{if} \hspace{18 pt}\quad\text{$i<j$,}\\
\text{Id},   \quad\quad\quad
 &\text{if}\hspace{17 pt} \quad   \text{$i=j$ or $i=j+1$,}\\
\delta_{i-1} \sigma_{j},  \quad
 &\text{    if} \hspace{16 pt}\quad \text{$i>j+1$}.
\end{cases}
\end{split}
\end{eqnarray}

A cocyclic module  is a cosimplicial module equipped with extra morphisms, $\tau :C^n\longrightarrow C^n$, called cocyclic maps such that the following extra relations hold.
\begin{eqnarray}
\begin{split}
&\tau\delta_{i}=\delta_{i-1} \tau, \hspace{43 pt} 1\le i\le  n+1,\\
&\tau \delta_{0} = \delta_{n+1}, \hspace{43 pt} \tau \sigma_{i} = \sigma _{i-1} \tau,  \hspace{33 pt} 1\le i\le n,\\ \label{rel2}
&\tau \sigma_{0} = \sigma_{n} \tau^2, \hspace{43 pt} \tau^{n+1} = \Id.
\end{split}
\end{eqnarray}
In a dual manner, one can define a cyclic module as a simplicial module with extra cyclic maps. More precisely:
 A cyclic module is given by data,  $ C=(C_{n},\delta_i, \sigma_{i}, \tau_{n})$, where $C_{n}$, $n\geq 0$ is a $\mathbb{C}$-module and  $\delta_{i}: C_{n}\longrightarrow C_{n-1},  \quad \sigma_{i}:C_n \ra C_{n+1},  \quad 0 \leq i \leq n$, and  $\tau:C_n \ra C_n$, are called faces, degeneracies and cyclic maps, respectively, are $\mathbb{C}$-module maps, satisfying the following  relations:
\begin{eqnarray}
\begin{split}
 &  \delta_{i} \delta_{j}=  \delta_{j-1} \delta_{i},  \hspace{28 pt}  \text{if} \quad i <j,\\
& \sigma_{i}  \sigma_{j} = \sigma_{j+1} \sigma_{i},  \hspace{25 pt} \text{if} \quad i \leq j,\\
&\delta_{i} \sigma_{j} = \label{rel11}
 \left\{\begin{matrix}
\sigma_{j-1} \delta_{i},   \hspace{10 pt}\text{if} \hspace{20 pt} i<j,\\
\Id,    \hspace{32 pt} \text{if}   \hspace{22 pt} i=j ~~ \text{or}~~ &i=j+1,\\
\sigma_{j} \delta_{i-1},  \hspace{8 pt}\text{if} \hspace{8 pt} i>j+1.
\end{matrix}\right.
\end{split}
\end{eqnarray}
and
\begin{align}
\begin{split}
&\delta_{i}\tau= \tau\delta_{i-1}, \hspace{43 pt} 1\le i\le n ,\\
& \delta_{0} \tau= \delta_{n} , \quad \sigma_{i} \tau= \tau  \sigma _{i-1},\hspace{23 pt} 1\le i\le n ,\\
&\sigma_{0} \tau = \tau^2  \sigma_{n} , \quad\tau^{n+1} = \Id. \\ \label{rel22}
\end{split}
\end{align}
\begin{remark}{\rm [Duality] Let $C=(C^{n}, d_i, s_i, t)$ be a cocyclic  module. We denote its cyclic dual by $\widetilde{C}$. It is shown in \cite{NCG} that with the following operators  $\widetilde{C}_n := C^{n}$ is a cyclic module.
\begin{align}
\begin{split}
&\d_i:=s_i: \widetilde{C}_n\longrightarrow \widetilde{C}_{n-1}, \quad 0\leq i\leq n-1,\quad \d_n:= \d_0 \tau_n,\\
&\sigma_i:=d_{i+1}:\widetilde{C}_n\longrightarrow \widetilde{C}_{n+1}, \quad 0\le i\le n-1,\\
&\tau:= t^{-1}.
\end{split}
\end{align}
Conversely, one obtains a cocyclic  module from a cyclic  module as follows. Let $C= (C_{n},\d_i, \sigma_i, \tau)$ be a cyclic module. Its cyclic dual is denoted by  $\breve{C}$  and defined by $\breve{C}^{n}= C_n $  with the following  cofaces, codegeneracies and cyclic operator.
\begin{align}
\begin{split}
& d_0: \tau_n \sigma_{n-1}, \quad d_i:= \sigma_{i-1}: \breve{C}^{n}\longrightarrow \breve{C}^{n+1}, \quad\quad 1\leq i \leq n,\\
&s_i:= \delta_i: \breve{C}^{n}\longrightarrow \breve{C}^{n-1}, \quad\quad 0\leq i \leq n-1, \\
& t:= \tau^{-1}.
\end{split}
\end{align}}
\end{remark}
%%%%%%%%%%%%%%%%%%%%%%%%%%%%%%%%%%%%%%%%%%%%%%%%%%%%%%%%%%%%%%%%%
\section{Hopf cyclic homology of $\times_{R}$-Hopf algebras}\label{S-2}
We refer the reader to \cite{b1,bwbook,sch,Scha:bia_nc} for details about bialgebroids and $\times$-Hopf algebras, however in Subsection \ref{SS-2-1} we   recall the basics of  $\times_R$-Hopf algebras. We review the stable anti Yetter-Drinfeld modules for $\times_R$-Hopf algebras in Subsection \ref{SS-2-2} from \cite{bs2}.  Finally in Subsection \ref{SS-2-3} we  bring the cyclic modules of module algebras and module corings over  $\times_R$-Hopf algebras with coefficients in stable anti Yetter-Drinfeld modules.
\subsection{ Preliminaries of $\times_{R}$-Hopf algebras}
\label{SS-2-1}
 Let $R$ and $\Kc$ be  algebras over  $\mathbb{C}$. In addition let the source and the target maps  be the $\Cb$-algebra homomorphisms $\mathfrak{s}:R\longrightarrow \Kc$ and $\mathfrak{t}:R^{op}\longrightarrow \Kc$, such that their
  ranges commute with one another. We equip  $\Kc$  with an $R$-bimodule structure  by $r_{1}kr_{2}=\mathfrak{s}(r_{1})\mathfrak{t}(r_{2})k$ for $r_1, r_2\in R$ and $k\in \Kc$.  Similarly  $\Kc\ot_R \Kc$  is endowed with the natural  $R$-bimodule structure, i.e,   $r_{1}(k_{1}\otimes_{R}k_{2})r_{2}:= \mathfrak{s}(r_{1})k_{1}\otimes_{R} \mathfrak{t}(r_{2})k_{2}.$ We assume that there are $R$-bimodule  maps called  the coproduct $\Delta:\Kc\longrightarrow \Kc\otimes_{R} \Kc$ and the couint $\varepsilon:\Kc\longrightarrow R$ via which   $\Kc$ is an $R$-coring. The data $(\Kc,\mathfrak{s},\mathfrak{t},\D,\ve)$ is called a left $R$-bialgebroid if  for $k_1,k_2\in \Kc$ and $r\in R$ the following identities hold

\begin{enumerate}
 \item [i)] $k^{(1)} \mathfrak{t}(r)\otimes_{R} k^{(2)} = k^{(1)} \otimes_{R} k^{(2)} \mathfrak{s}(r),$
\item[ii)] $\Delta(1_{\Kc})=1_{\Kc}\ot_R 1_{\Kc}, \quad \text{ and} \quad  \Delta(k_1k_2)=k_1^{(1)} k_2^{(1)}\otimes_R k_1^{(2)} k_2^{(2)}$.
\item[iii)] $\varepsilon(1_{\Kc})=1_{R} \quad \text{and} \quad \varepsilon(k_1k_2)=\varepsilon(k_1\mathfrak{s}(\varepsilon(k_2))).$
\end{enumerate}
%\begin{example}\label{ex1} {\rm Let $B$ be an algebra over $\mathbb{C}$. The $B^{e}=B\ot B^{op}$ is a left $\times_B$-bialgebroid where the source and target maps can be %defined by
%$$s: B\longrightarrow B^{e}, \quad b\longmapsto b\ot 1; \quad t: B^{op}\longrightarrow b^{e}, \quad b^{o}\longmapsto 1\ot b^{o},$$  comultiplication by
%$$\Delta:B^{e}\longrightarrow B^{e}\ot_B B^{e}, \quad b\ot b^{o}\longmapsto (b\ot 1)\ot_B (1\ot b^{o}),$$ and counit by
%$$ \varepsilon:B^{e}\longrightarrow B, \quad \varepsilon(b\ot b')=bb'. $$}
%\end{example}

 For a left bialgebroid $\Kc$ over an algebra $R$, the following identities hold for  $r, r_1, r_2, r_3, r_4\in R$ and $k, k_1, k_2\in \Kc$ \cite[4.31]{bwbook}.

\begin{enumerate}
\item[a)] $\varepsilon(\mathfrak{s}(r))= \varepsilon(\mathfrak{t}(r)) = r.$
\item[b)] $\mathfrak{s}(\varepsilon(k^{(1)}))k^{(2)}= \mathfrak{t}(\varepsilon(k^{(2)}))k^{(1)}=k.$
\item[c)] $\varepsilon(k_1k_2)=\varepsilon(k_1\mathfrak{t}(\varepsilon(k_2))).$
\item[d)] $\Delta(\mathfrak{s}(r))=\mathfrak{s}(r)\otimes_{R} 1_\Kc \quad \text{and} \quad \Delta(\mathfrak{t}(r))=1_\Kc\otimes_{R} \mathfrak{t}(r).$
\item[e)]$\Delta(\mathfrak{s}(r_1)\mathfrak{t}(r_2)k\mathfrak{s}(r_3)t(r_4))=\mathfrak{s}(r_1)k^{(1)}\mathfrak{s}(r_3)\otimes_{R}\mathfrak{t}(r_2)k^{(2)}\mathfrak{t}(r_4).$
\end{enumerate}
\medskip

Similarly  a right $R$-bialgebroid $\Bc$ is an $R$-coring $(\Bc,\mathfrak{s},\mathfrak{t},\D,\ve)$, where $\mathfrak{s},\mathfrak{t}:R\ra \Bc$ are algebra and anti-algebra homomorphisms  with commuting ranges, and $\Bc$ is an  $R$-bimodule with the right multiplications by $\mathfrak{s}$ and $\mathfrak{t}$, such that the following conditions are satisfied for $b,b_1,b_2\in \Bc$ and $r\in R$.
\begin{enumerate}
\item[i)]  $b^{(1)}\ot_R \mathfrak{t}(r)b^{(2)}=\mathfrak{s}(r)b^{(1)}\ot_R b^{(2)},$
\item[ii)] $ \D(1_B)=1_B\ot_R1_B \quad \text { and} \quad \D(b_1b_2)=b_1^{(1)}b_2^{(1)}\ot_R b_1^{(2)}b_2^{(2)}, $
\item[iii)] $ \varepsilon(1_B)=1_R \quad \text{and} \quad  \varepsilon(b_1b_2)=\varepsilon(\mathfrak{s}(\ve(b_1))b_2).$
\end{enumerate}

%\begin{example}\cite[page 11, 3.2.3]{b1}\label{ex2} {\rm  Let $B$ be an algebra over  $\mathbb{C}$. The $B^{e}=B\ot B^{op}$ is a right $\times_B$-bialgebroid where the source and target maps are defined by
%$$s: B\longrightarrow B^{e}, \quad b\longmapsto b\ot 1; \quad t: B^{op}\longrightarrow b^{e}, \quad b^{o}\longmapsto 1\ot b^{o},$$  comultiplication by
%$$\Delta:B^{e}\longrightarrow B^{e}\ot_B B^{e}, \quad b\ot b^{o}\longmapsto (1\ot b^{o})\ot_B (b\ot 1),$$ and counit by
%$$ \varepsilon:B^{e}\longrightarrow B, \quad \varepsilon(b\ot b')=bb'. $$}
%\end{example}
 For a right bialgebroid $\mathcal{B}$ over an algebra $R$, the following identities hold for $r, r_{1}, r_2, r_3, r_4 \in R$ and $b, b_1, b_2\in B$.

\begin{enumerate}
\item[a)] $\varepsilon(\mathfrak{s}(r))= \varepsilon(\mathfrak{t}(r)) = r.$
\item[b)] $b=b^{(1)}\mathfrak{s}(\varepsilon(b^{(2)}))=b^{(2)}\mathfrak{t}(\varepsilon(b^{(1)})).$
\item[c)] $\varepsilon(b_1b_2)=\varepsilon(\mathfrak{t}(\varepsilon(b_1))b_2).$
\item[d)] $\Delta(\mathfrak{s}(r))=1_B \otimes_{R} \mathfrak{s}(r) \quad \text{ and}  \quad \Delta(\mathfrak{t}(r))=\mathfrak{t}(r)\otimes_{R} 1_B.$
\item[e)] $\Delta(\mathfrak{s}(r_1)\mathfrak{t}(r_2)b\mathfrak{s}(r_3)\mathfrak{t}(r_4))=\mathfrak{t}(r_2)b^{(1)}\mathfrak{t}(r_4)\otimes_{R}\mathfrak{s}(r_1)b^{(2)}\mathfrak{s}(r_3).$
\end{enumerate}

 A left $R$-bialgebroid $(\Kc,\mathfrak{s},\mathfrak{t},\Delta,\ve)$ is said to be a left $\times_{R}$-Hopf algebra provided that the map
\begin{align}\label{aa}
\nu: \Kc\otimes_{R^{op}}\Kc\longrightarrow \Kc\otimes_{R}\Kc, \qquad  k\otimes_{R^{op}}k'\longmapsto k^{(1)}\otimes_{R}k^{(2)}k',
\end{align}
is bijective. Here in the  left hand side of \eqref{aa} the  $R^{op}$-module structures are given by right and left multiplication by $\mathfrak{t}(r)$ for $r\in R$, however in the right hand side the $R$-module structures are given by  the original actions  of $R$ on $\Kc$.
One notes that $\nu$ and $\nu^{-1}$ are right $\mathcal{K}$-linear. The role of antipode in $\times_R$-Hopf algebras is played by the following map.
\begin{align}
\begin{split}
\Kc\ra \Kc\ot_{R^{op}}\Kc, \quad k\mapsto \nu^{-1}(k\ot_{R} 1_{\mathcal{K}}).
\end{split}
\end{align}

 For our convenience we use $k^-\ot_{R^{op}} k^+ $ for $\nu^{-1}(k\otimes_{R}1_{\mathcal{K}})$.

\begin{example}\label{envelop}{\rm Let $R$ be an algebra over $\mathbb{C}$. The simplest left $\times_R$-Hopf algebra is $\mathcal{K}=R\ot R^{op}$  with the source and  target maps defined by
$$\mathfrak{s}: R\longrightarrow \mathcal{K}, \quad r\longmapsto r\ot 1; \quad \mathfrak{t}: R^{op}\longrightarrow \mathcal{K}, \quad r\longmapsto 1\ot r,$$  comultiplication defined by
$$\Delta:\mathcal{K}\longrightarrow \mathcal{K}\ot_R \mathcal{K}, \quad r_1\ot r_2\longmapsto (r_1\ot 1)\ot_R (1\ot r_2),$$  counit given by
$$ \varepsilon:\mathcal{K}\longrightarrow R, \quad \varepsilon(r_1\ot r_2)=r_1r_2,  $$ and
$$\nu((r_1\ot r_2)\ot(r_3\ot r_4))= r_1\ot 1\ot r_3\ot r_4r_2,$$
$$\nu^{-1}((r_1\ot r_2)\ot(r_3\ot r_4))=r_1\ot 1\ot r_2r_3\ot r_4, $$ where $r,r_1, r_2, r_3, r_4 \in R.$
}
\end{example}

The following properties are a symmetrical version of \cite[Proposition 3.7]{sch}. We refer the reader to \cite[Lemma 2.14]{bs2} for more properties of a right$\times$-Hopf algebras.

\begin{proposition}\label{lem} Let $R$ be an algebra over $\mathbb{C}$ and $\mathcal{K}$ a left $\times_{R}$-Hopf algebra. The following identities hold for  $k,k'\in \mathcal{K}$ and $r\in R$.
\begin{enumerate}
\item[\rm i)] ${k^{-}}^{(1)}\otimes_{R}{k^{-}}^{(2)}k^{+}=k\otimes_{R}1_{\mathcal{K}}.$
 \item[\rm ii)] ${k^{(1)}}^{-}\otimes_{R^{op}}{k^{(1)}}^{+}k^{(2)}=k\otimes_{R^{op}}1_{\mathcal{K}}.$
 \item[\rm iii)] $(kk')^{-}\otimes_{R^{op}}(kk')^{+}= k^{-}k'^{-}\otimes_{R^{op}}k'^{+}k^{+}.$
 \item[\rm iv)] $1_{\mathcal{K}}^{-}\otimes_{R^{op}}1_{\mathcal{K}}^{+}=1_{\mathcal{K}}\otimes_{R^{op}}1_{\mathcal{K}}.$
 \item[\rm v)] $ {k^{-}}^{(1)}\otimes_{R}{k^{-}}^{(2)}\otimes_{R^{op}}k^{+}= k^{(1)}\otimes_{R}{k^{(2)}}^{-}\otimes_{R^{op}}{k^{(2)}}^{+}.$
\item[\rm vi)] $k^{-}\ot_{R^{op}} k^{{+}(1)}\ot_{R} k^{{+}(2)}=    k^{--}\ot_{R^{op}} k^{+} \ot_R k^{-+}.$
 \item[\rm vii)] $k= k^{-}\mathfrak{t}(\varepsilon(k^{+})).$
 \item[\rm viii)] $\mathfrak{s}(\varepsilon(k))=k^{-}k^{+}$.
 \item[\rm ix)] $k^{-}\otimes_{R^{op}}k^{+}\mathfrak{t}(r)=\mathfrak{t}(r)k^{-}\otimes_{R^{op}}k^{+}.$
\end{enumerate}
\end{proposition}

 Let $R$ be an  algebra over $\mathbb{C}$. A right $R$-bialgebroid $\mathcal{B}=(B,\mathfrak{s},\mathfrak{t},\Delta,\varepsilon)$, is said to be a right $\times_{R}$-Hopf algebra provided that the map
\begin{align}\label{map}
\nu: \mathcal{B}\otimes_{R^{op}}\mathcal{B}\longrightarrow \mathcal{B}\otimes_{R}\mathcal{B}, \qquad  b\otimes_{R^{op}}b'\longmapsto bb'^{(1)}\otimes_{R}b'^{(2)},
\end{align}
is bijective. In the domain of the map \eqref{map}, $R^{op}$-right  and left module structures of $\mathcal{B}$ are given by right and left multiplication by $\mathfrak{t}(r)$ respectively, for $r\in R$. In the codomain of the map \eqref{map}, $R$-module structures are given by right multiplication by $\mathfrak{s}(r)$ and $\mathfrak{t}(r)$.

 One notes that $\nu$ and $\nu^{-1}$ are   left $\mathcal{B}$-linear maps. We denote the image of $\nu^{-1}$ by the following index notation, $$\nu^{-1}(1\otimes_{R^{op}}b)=b^{-}\otimes_{R}b^{+}.$$

\begin{example}{\rm Let $R$ be an algebra over $\mathbb{C}$. The $\mathcal{B}=R\ot R^{op}$ is a right $\times_R$-Hopf algebra where the source and target maps are defined by
$$\mathfrak{s}: R\longrightarrow \mathcal{B}, \quad r\longmapsto r\ot 1; \quad \mathfrak{t}: R^{op}\longrightarrow \mathcal{B}, \quad r\longmapsto 1\ot r,$$  comultiplication by
$$\Delta:\mathcal{B}\longrightarrow \mathcal{B}\ot_B \mathcal{B}, \quad r_1\ot r_2\longmapsto (1\ot r_2)\ot_B (r_1\ot 1),$$  counit by
$$ \varepsilon:\mathcal{B}^{e}\longrightarrow B, \quad \varepsilon(r_1\ot r_2)=r_2r_1, $$ and
$$\nu((r_1\ot r_2)\ot (r_3\ot r_4))= r_1\ot r_4r_2\ot r_3\ot 1,$$
$$\nu^{-1}((r_1\ot r_2)\ot (r_3\ot r_4))= r_1 r_4\ot r_2\ot r_3\ot 1.$$

}
\end{example}

%%%%%%%%%%%%%%%%%%%%%%%%%%%%%%%%%%%%%%%%%%%%%%%%%%%%%%%%%%%%%%%%%%%%%%%%%%%%%%%%%%%%%%%%%%%%%%%%%%%%%%%%%%%%%%
\subsection{Stable anti Yetter-Drinfeld-modules}
\label{SS-2-2}
In this subsection for the reader's convenience we briefly recall the definitions of  module, comodule and stable anti Yetter-Drinfeld-module for  $\times_{R}$-Hopf algebras. Also we present an example of a stable anti Yetter-Drinfeld (SAYD) module for the $\times_R$-Hopf algebra $R\ot R^{op}$ which plays an important role in the sequel  section.

A right module over a  $\times_{R}$-Hopf algebra $\mathcal{K}$ is  a right $\mathcal{K}$-module $M$, which is naturally an  $R$-bimodule  by $rm=m\triangleleft \mathfrak{t}(r),$  and $mr=m \triangleleft \mathfrak{s}(r)$. Similarly, by a  left module over a  right $\times_{R}$-Hopf algebra $\mathcal{B}$ we mean a left $\Bc$-module $M$, which  is naturally an $R$-bimodule by $rm=\mathfrak{s}(r)\triangleright m,$  and  $ mr=\mathfrak{t}(r)\triangleright m$.

A left (right) comodule  over an  $R$-coring $\Cc$  is defined by  a left(right) $R$-module $M$, together with a left(right) $R$-module coaction map,
\begin{equation}
M\longrightarrow \Cc\otimes_{R}M,\quad m\longmapsto m\ns{-1}\otimes_{R}m\ns{0},
 \end{equation}
 ( $M\longrightarrow M\ot_R \Cc,$ $m\longmapsto m\ns{0}\ot_R m\ns{1}$) satisfying coassociativity and counitality axioms.

By a comodule over a  $\times$-Hopf algebra we mean a comodule over the underlying coring. One notes that a left comodule $M$ of a left bialgebroid $\mathcal{K}$ can be equipped with an $R$-bimodule structure by introducing a right $R$-action
\begin{align}\label{act}
m r:=\varepsilon(m\ns{-1}\mathfrak{s}(r))\triangleright m\ns{0},
\end{align}
for $r\in R$ and $m\in M$. It is checked in \cite{bs2} that the left $\mathcal{K}$-coaction on left comodule $M$ is an $R$-bimodule map. That is  for $r,r'\in R$ and $m\in M$,
\begin{align}
(rmr')\ns{-1}\otimes_{R} (rmr')\ns{0}= \mathfrak{s}(r)m\ns{-1}\mathfrak{s}(r')\otimes_{R} m\ns{0}.
\end{align}
Furthermore,  for any $m\in M$ and $r\in R$,
\begin{align}\label{abc}
m\ns{-1}\otimes_{R} m\ns{0}r= m\ns{-1}\mathfrak{t}(r)\otimes_{R}m\ns{0}.
\end{align}

Here we recall the definition of  stable anti-Yetter-Drinfeld modules over $\times$-Hopf algebras from \cite{bs2}.
  Let $\Kc$ be a left $\times_R$-Hopf algebra and $M$ be a right $\mathcal{K}$-module and a left $\mathcal{K}$-comodule.  We say $M$ is a right-left anti Yetter-Drinfeld module provided that the following conditions hold.

\begin{itemize}

\item [i)] The $R$-bimodule structures on $M$, underlying its module and comodule structures, coincide. That is, for $m\in M$ and $r\in R$,
$$mr= m\triangleleft \mathfrak{s}(r), \qquad \text{and} \qquad rm= m \triangleleft \mathfrak{t}(r),$$ where $rm$ denotes the left $R$-action on the left $\mathcal{K}$-comodule $M$ and $mr$ is the canonical right action  defined in\eqref{act}, i.e. $mr= \varepsilon(m\ns{-1}\mathfrak{s}(r))m\ns{0}.$

\item[ii)] For $k\in\mathcal{K}$ and $m\in M$,
\begin{equation}\label{b}
(m\triangleleft k)\ns{-1}\otimes(m \triangleleft k)\ns{0}= {k^{(2)}}^{+} m\ns{-1}k^{(1)}\otimes_{R} m\ns{0}\triangleleft {k^{(2)}}^{-} .
\end{equation}
\end{itemize}
The anti Yetter-Drinfeld module $M$ is said to be stable if in addition, for any $m\in M$, $m\ns{0}m\ns{-1}=m.$

Similarly one defines  a   left-right AYD module  over a right $\times_R$-Hopf algebra  $\Bc$ by a left $\mathcal{B}$-module and  right $\mathcal{B}$-comodule $M$ satisfying
\begin{itemize}
\item[i)] $mr= \mathfrak{t}(r)\triangleright m, \qquad \text{and} \qquad rm= \mathfrak{s}(r)\triangleright m,$ where $mr$ denotes the right $R$-action on the right $\mathcal{B}$-comodule $M$ and $rm$ is the canonical left $R$-action $rm=m\ns{0}. \varepsilon(\mathfrak{s}(r)m\ns{1})$.
\item[ii)] For $b\in\mathcal{B}$ and $m\in M$,
\begin{align}\label{bb}
(b\triangleright m)\ns{0}\otimes(b\triangleright m)\ns{1}= {b^{(1)}}^{+}\triangleright  m\ns{0}\otimes_{R} b^{(2)} m\ns{1} {b^{(1)}}^{-} .
\end{align}
\end{itemize}

The anti Yetter-Drinfeld module $M$ is said to be stable if in addition, for any $m\in M$, we have $m\ns{1}m\ns{0}=m.$

\begin{remark}{\rm Left or right $\times_{R}$-Hopf algebras extend the notion of   Hopf algebras. In fact if $\mathcal{K}$ is a bialgebra over a commutative ring $R$, with coproduct $k\longmapsto k^{(1)}\otimes_{R}k^{(2)},$ then the bijectivity of the map $\nu$ is equivalent to the fact that $\mathcal{K}$ is a Hopf algebra. In this case the inverse of the map $\nu$ can be defined as
$$\nu^{-1}(k\otimes 1)= k^{-}\otimes_{R^{op}} k^{+}:= k^{(1)}\otimes_{R^{op}}S(k^{(2)}),$$ where $S$ denotes the antipode of the Hopf algebra $\mathcal{K}$. Therefore the condition \eqref{b} for a left $\times_{R}$-Hopf algebra is equivalent to the following relation
$$(mk)\ns{-1}\otimes (mk)\ns{0}=S(k^{(3)})m\ns{-1}k^{(1)}\otimes m\ns{0}k^{(2)}.$$ }
\end{remark}
Here we recall the definition of a character for $\times_R$-Hopf algebras.
 A map  $\d$ is called a  right character \cite[Lemma 2.5]{b1}, for the  $\times_R$-Hopf algebra $\mathcal{K}$ if it satisfies the following conditions:
\begin{align}\label{char-1}
& \delta(k \mathfrak{s}(r))= \delta(k)r, \quad \text{for} \quad k\in \mathcal{K} \quad \text{and} \quad r\in R,\\\label{char-2}
& \delta(k_1k_2)=\delta(\mathfrak{s}(\delta(k_1))k_2), \quad \text{ for} \quad k_1,k_2\in \mathcal{K},\\\label{char-3}
& \delta(1_{\mathcal{K}})=1_R.
\end{align}
As an example, for any right $\times_R$-Hopf algebra the counit $\varepsilon$ is a right character. The following example is similar to \cite[Example 2.18]{bs2} for right $\times_R$-Hopf algebras.

\begin{example}\label{identity}{\rm
Let $\mathcal{K}$ be a left $\times_R$-Hopf algebra, $\sigma \in \mathcal{K}$  a group-like element and the map $\d:\mathcal{K}\longrightarrow R$ a right character.
The following action and coaction,
\begin{align}\label{RSAYD}
r \triangleleft k= \d(\mathfrak{s}(r)k), \quad \text{and} \quad r\longmapsto \mathfrak{s}(r)\sigma \ot 1
\end{align}

define a right $\mathcal{K}$-module and left $\mathcal{K}$-comodule structure  on $R$, respectively. These action and coaction amount  to a right-left anti Yetter-Drinfeld module on $R$ if and only if, for all $r\in R$ and $k\in \mathcal{K}$,
\begin{align}
\mathfrak{s}(\d(k))\sigma=\mathfrak{t}(\d(k^{(2)-}))k^{(2)+}\sigma k^{(1)}, \quad \text{and} \quad \varepsilon(\sigma \mathfrak{s}(r))=\d(\mathfrak{s}(r)).
\end{align}
The anti Yetter-Drinfeld module $R$ is stable if in addition $\d(\mathfrak{s}(r)\sigma)=r$, for all $r\in R$.  We specialize this example to the left $\times_R$-Hopf algebra $\mathcal{K}=R\ot R^{op}$ explained in the Example \ref{envelop}.
}
\end{example}

To this end,  we like to know all group-like elements and right characters of  $\Kc$.  One can easily characterize all homogenous group-like elements of $\Kc$. An element $  x\ot y \in \Kc$ is a group-like element if and only if,  $ xy=yx=1$, because
\begin{align*}
&\Delta( x\ot y)= (x\ot y)\ot_R  (x\ot y)= (x\ot y)\ot_R (x\ot1)(1\ot y)=\\
& (1\ot x)(x\ot y)\ot_R (1\ot y)=(x\ot yx) \ot_R (1\ot y).
\end{align*}
On the other hand,  one sees $ \Delta(x\ot y)=(x\ot 1) \ot_R (1\ot y).$ Therefore $yx=1$. Since $x\ot y$ is a group like element,  we have
$1=\varepsilon( x\ot y)= xy.$ Conversely one easily sees that if $ xy=yx =1$ then $ x\ot y$ is a group-like element for $\mathcal{K}$.

        We claim that a map  $\delta:\mathcal{K}\longrightarrow R$  is a right character  if and only if   $\delta(s\ot r)= \theta(s)t$ for some algebra map $\theta: R\ra R$. Indeed,  let $\d$ be a right character on $\mathcal{K}$,  then
  $$\delta(s\ot r)=\delta((1\ot r)(s\ot 1))=\delta((1\ot r)\mathfrak{s}(s))=\d(1\ot r)s.$$
  We define $\t$ by $\theta(r):=\d(1\ot r)$. It is obvious that $\t$ is  unital. We show that $\theta$ is an algebra map,
\begin{align*}
\theta(rr') &=\d(1\ot rr')=\d((1\ot r')(1\ot r))=\d(\mathfrak{s}(\d((1\ot r'))(1\ot r))\\
&=\d((1\ot r)\mathfrak{s}(\d((1\ot r'))=\d(1\ot r)\d(r'\ot 1)=\theta(r)\theta(r').
\end{align*}
 Conversely, let $\t$ be an algebra endomorphism of $R$, we consider the
 map $\d:\mathcal{K}\longrightarrow R, \quad s\ot r \longmapsto \theta(r)s$. We
 show that $\d$ satisfies \eqref{char-1}. Indeed, $\d((s\ot r)\mathfrak{s}(r'))=\d(sr'\ot r )=\theta(r)sr'=\d(s\ot r)r'.$
The satisfaction of \eqref{char-2} by $\d$ is as follows.
\begin{align*}
&\d((s\ot r)(s'\ot r'))=\d(ss'\ot r'r)=\theta(r'r)ss'=\theta(r')\theta(r)ss'=\d(\theta(r)ss'\ot r'))=\\
&\d((\theta(r)s\ot 1)(s' \ot r'))=\d((\mathfrak{s}(\theta(r)s)(s',r')=\d(\mathfrak{s}(\d(s,r))(s'\ot r'))).
\end{align*}
Also since $\d(1_R\ot 1_R)=\t(1_R)1_R=1_R$, the map $\d$ satisfies \eqref{char-3}. Therefore $\d$ is a right character.

\begin{proposition}\label{nice}
 Let $\mathcal{K}$ be  the left $\times_{R}$-Hopf algebra $R\ot R^{op}$,  $x\ot x^{-1}\in \mathcal{K}$ a homogenous group like element and  $\delta$ a right character on $\mathcal{K}$. Then the  action and  coaction defined in \eqref{RSAYD}
 amount  to a right-left SAYD module on $R$ if and only if  $x$ belongs the center of the algebra $R$ and $\theta=Id$.
\end{proposition}

\begin{proof} Let  $k\ot t\in \mathcal{K}$ and $r\in R$. By the characterization of all right characters on $\mathcal{K}$ the action and coaction defined in \eqref{RSAYD} reduce to
$$ r \triangleleft (k\ot t)= \theta(t)rk, \quad \text{and} \quad r\longmapsto (rx\ot x^{-1})\ot 1.$$ Now we show   $x$  is in  the center of the algebra $R$. For this we use the stability condition,
$$r= \d(\mathfrak{s}(r)(x\ot x^{-1}))=\d((r\ot 1)(x\ot x^{-1}))= \d(rx\ot x^{-1}), \quad \forall r\in R.$$ In the previous equality let $r=x^{-1}$. We have
$$y=\d(1\ot x^{-1})= \theta(x^{-1}).$$ On the other hand,
$$r= \d(\mathfrak{s}(r)(x\ot x^{-1}))= \d(rx\ot x^{-1})= \theta(x^{-1})rx.$$ By multiplying the both sides by $x^{-1}$, we obtain $\theta(x^{-1})r=rx^{-1}$. Therefore $x^{-1}r=rx^{-1}$ for all $r\in R$. Thus $x$ is also in the center of the algebra $R$.

 We shall show the AYD condition implies that the map $\theta$ is identity.
 The AYD condition is equivalent to
$$\mathfrak{s}(\d(k))\sigma=t(\d({k^{(2)}}^{-})){k^{(2)}}^{+}\sigma k^{(1)}.$$ Let $k=g\ot h$ and $\sigma=x\ot x^{-1}$. We have
$$\theta(h)gx\ot x^{-1}= (h\ot 1)(x\ot x^{-1})(g\ot 1)=hxg\ot x^{-1},$$ and therefore $\theta(h)gx=hxg$. Since $x$ is in the center of $R$,  $\theta(h)gx=hgx$. By multiplying both sides by $x^{-1}$ from right, we obtain $\theta(h)g=hg$ and for $g=1$ we have $\theta(h)=h$ for $\forall h\in R$.
Thus we obtain the claimed action and coaction.\\
Conversely, it is easy to check that if $x$ belongs the center of $R$ and $\theta=Id$, then the following  $R\ot R^{op}$-action and  coaction
\begin{equation}
 r \triangleleft (k\ot t)= trk, \quad \text{and} \quad r\longmapsto (rx\ot x^{-1})\ot 1, \quad  k,t,r\in R,
 \end{equation}
 define a right-left SAYD module $R$ over the left $\times_R$-Hopf algebra $ R\ot R^{op}$.

\end{proof}

%%%%%%%%%%%%%%%%%%%%%%%%%%%%%%%%%%%%%%%%%%%%%%%%%%%%%%%%%%%%%%%%%%%%%%%%%%%%%%%%%%%%%%%%%%%%%%%%%%%%%%%%%%%%%%%%%%%%%%%%%%%%
\subsection{Cyclic homology of module algebras and module corings}
\label{SS-2-3}
In this subsection we recall from \cite{bs2} the cyclic homology  of  module algebras and module corings with coefficients in  SAYD modules under the symmetry of  $\times_{R}$-Hopf algebras . In the case of module algebras,  first we introduce a cocyclic module for left $\times_{R}$-Hopf algebras and then we dualize it, as explained in the Section \ref{dual}  to find a cyclic module.

 Let $R$ be an algebra over   $\mathbb{C}$ and $\mathcal{K}$ be a left $\times_{R}$-Hopf algebra. A left $\mathcal{K}$-module algebra $A$ \cite[Definition 2.3]{bs2} is a $\mathbb{C}$-algebra and a left $\mathcal{K}$-module such that for all $k\in \mathcal{K}$ , $a,a'\in A$ and $r\in R$, the following identities hold,
\begin{align}\label{modulealgebra-1}
& k\triangleright 1_{A}=\mathfrak{s}(\varepsilon(k))\triangleright 1_{A},\\ \label{modulealgebra-2}
 &k\triangleright (a a')= (k^{(1)}\triangleright a)(k^{(2)}\vartriangleright a'),\\\label{modulealgebra-3}
 &(\mathfrak{t}(r)\triangleright a)a'=a(\mathfrak{s}(r)\triangleright a'),  \text{(multiplication is $R$-balanced.)}
\end{align}

\begin{proposition} Let $S$ be an algebra over $\mathbb{C}$, $\mathcal{K}$ a left $\times_{S}$-Hopf algebra, $T$ a left $\mathcal{K}$-module algebra,  and $M$ a right-left SAYD module over $\mathcal{K}$.  Let
\begin{equation}
_\Kc{C^{n}}(T, M)= M\otimes_{\mathcal{K}}T^{\otimes_{S}(n+1)}.\end{equation}
Then the  following cofaces, codegeneracies and cyclic map define a cocyclic module structure on $_\Kc{C^\ast}(T,M)$.
\begin{eqnarray}
&& d_{0}(m\otimes_{\mathcal{K}} \td t)= m\otimes_{\mathcal{K}}1\otimes_{S} t_{0}\otimes_{S}\cdots \otimes_{S}t_{n},\nonumber  \\
&& d_{i}(m\otimes_{\mathcal{K}} \td t)=m\otimes_{\mathcal{K}}t_{0}\otimes_{S}\cdots\otimes_{S}t_{i}\otimes_{S}1\otimes_{S}t_{i+1}\otimes\cdots\otimes_{S}t_{n}, \quad 1\leq i\leq n,  \nonumber \\
&& d_{n}(m\otimes_{\mathcal{K}} \td t)=m\ns{0}\otimes_{\mathcal{K}}t_{0}\otimes_{S}\cdots\otimes_{S}t_{n}\ot_S m\ns{-1}\triangleright 1_{T},   \nonumber \\
&& s_{i}(m\otimes_{\mathcal{K}}\td t)= m\otimes_{\mathcal{K}}t_{0}\otimes_{S}\cdots \otimes_{S}t_{i}t_{i+1}\otimes_{S}\cdots t_{n}, \quad 0\leq i\leq n-1,\nonumber\\
&& \mathrm{t}_n(m\otimes_{\mathcal{K}} \td t)= m\ns{0}\otimes_{\mathcal{K}}t_{1}\otimes_{S}\cdots \otimes_{S}t_{n}\otimes_{S} m\ns{-1}\triangleright t_{0}. \nonumber
\end{eqnarray}
Here $\td t$ stands for $t_0\ot_S\cdots\ot_S t_n$.
\end{proposition}
\begin{proof} Since $T$ is an $\mathcal{K}$-module algebra, by $S$-balanced property of $T$, we have $k\triangleright 1_{T}=\mathfrak{s}(\varepsilon(k))1_{T}=\mathfrak{t}(\varepsilon(k))1_{T}$. Therefore,
  cofaces are well-defined. Here we consider $T^{\otimes_{S}(n+1)}$ as a $\mathcal{K}$-module with diagonal action. By \eqref{modulealgebra-2}  the codegeneracies are well-defined. We show that the cocyclic map  is well-defined:
\begin{eqnarray}
&& \mathrm{t}_{n}(m\triangleleft k\otimes_{\mathcal{K}} t_{0}\otimes_{S}\cdots\otimes_{S}t_{n}) \nonumber \\
&=& (m\triangleleft k)\ns{0}\otimes_{\mathcal{K}} t_{1}\otimes_{S}\cdots \otimes_{S}t_{n}\otimes_{S} ((m\triangleleft k)\ns{-1})\triangleright t_0 \nonumber \\
&=& m\ns{0}\triangleleft {k^{(2)}}^{-}\otimes_{\mathcal{K}} t_{1}\otimes_{S}\cdots \otimes_{S}t_{n}\otimes_{S}{k^{(2)}}^{+}m\ns{0}k^{(1)} \triangleright t_0\nonumber \\
&=& m\ns{0}\triangleleft {k^{-}}^{(2)}\otimes_{\mathcal{K}} t_{1}\otimes_{S}\cdots \otimes_{S}t_{n}\otimes_{S} k^{+}m\ns{-1}{k^{-}}^{(1)}\triangleright t_0 \nonumber \\
&=& m\ns{0}\otimes_{\mathcal{K}} {k^{-}}^{(2)}\triangleright( t_{1}\otimes_{S}\cdots \otimes_{S}t_{n}\otimes_{S} k^{+}m\ns{-1}{k^{-}}^{(1)}\triangleright t_0) \nonumber \\
&=& m\ns{0}\otimes_{\mathcal{K}} {k^{-}}^{(2)(1)}\triangleright t_{1}\otimes_{S}\cdots \otimes_{S}{k^{-}}^{(2)(n)}\triangleright t_{n}\otimes_{S}\nonumber\\
&&{k^{-}}^{(2)(n+1)}k^{+}m\ns{-1}{k^{-}}^{(1)}\triangleright t_0 \nonumber \\
&=& m\ns{0}\otimes_{\mathcal{K}} {k^{-}}^{(2)}\triangleright t_{1}\otimes_{S}\cdots \otimes_{S}{k^{-}}^{(n+1)}\triangleright t_{n}\otimes_{S}{k^{-}}^{(n+2)}k^{+} m\ns{-1}{k^{-}}^{(1)}\triangleright t_0 \nonumber \\
&=& m\ns{0}\otimes_{\mathcal{K}}k^{(2)}\triangleright t_{1}\otimes_{S}\cdots \otimes_{S}k^{(n+1)}\triangleright t_{n}\otimes_{S}m\ns{-1}k^{(1)}\triangleright t_{0} \nonumber \\
&=& \mathrm{t}_{n}(m\otimes_{\mathcal{K}} k^{(1)}\triangleright t_{0}\otimes_{S}\cdots \otimes_{S}k^{(n+1)}\triangleright t_{n})\nonumber \\
&=& \mathrm{t}_{n}(m \otimes_{\mathcal{K}}k\triangleright( t_{0}\otimes_{S}\cdots\otimes_{S}t_{n})). \nonumber
\end{eqnarray}

In the above, we use the SAYD condition \eqref{b} in the second equality, Proposition \ref{lem}(v) in the third equality,  the diagonal action of $\mathcal{K}$ on $T^{\otimes_S (n+1)}$ in the fifth equality and  Proposition \ref{lem}(i) in the seventh equality. Now we check some of the cocyclicity conditions. To see $\mathrm{t}_{n}^{n+1}= Id$, we have
\begin{eqnarray}
&& \mathrm{t}_{n}^{n+1}(m\ot_{\mathcal{K}}t_0\ot_S \cdots \ot_S t_{n+1}) \nonumber \\
 &=& m\ns{0}\ot_{\mathcal{K}}m\ns{-n-1}t_0\ot_S\cdots \ot_S m\ns{-1}t_n \nonumber \\
 &=& m\ns{0}\ot_{\mathcal{K}}m\ns{-1}\triangleright (t_0\ot_S\cdots \ot_S t_n )\nonumber \\
 &=& m\ns{0}m\ns{-1}\ot_{\mathcal{K}}t_0\ot_S\cdots \ot_S t_n \nonumber \\
 &=& m \ot_{\mathcal{K}}t_0\ot_S\cdots \ot_S t_n .\nonumber \\
\end{eqnarray}
 In the above, we use the left diagonal action of $\mathcal{K}$ on $T^{\ot_{S}(n+1)}$ in the second equality and the stability condition in the last equality.

Now we show $\mathrm{t}_{n}s_{0}= s_{n} \mathrm{t}_{n+1}^{2}$.
\begin{eqnarray}
&&  s_{n} \mathrm{t}_{n+1}^{2}(m\ot_{\mathcal{K}}t_0\ot_S \cdots \ot_S t_{n+1}) \nonumber \\
&=&s_n \mathrm{t}_{n+1}(m\ns{0}\ot_{\mathcal{K}}t_1\ot_S\cdots \ot_S t_{n+1}\ot_S m\ns{-1}\triangleright t_0)\nonumber \\
&=&s_n (m\ns{0}\ns{0}\ot_{\mathcal{K}}t_2\ot_S\cdots \ot_S t_{n+1}\ot_S m\ns{-1}\triangleright t_0 \ot_S m\ns{0}\ns{-1}\triangleright t_1)\nonumber \\
&=&s_n (m\ns{0}\ot_{\mathcal{K}}t_2\ot_S\cdots \ot_S t_{n+1}\ot_S {m\ns{-1}}^{(1)}\triangleright t_0 \ot_S {m\ns{-1}}^{(2)}\triangleright t_1)\nonumber \\
&=&m\ns{0}\ot_{\mathcal{K}}t_2\ot_S\cdots \ot_S t_{n+1}\ot_S ({m\ns{-1}}^{(1)}\triangleright t_0 )( {m\ns{-1}}^{(2)}\triangleright t_1)\nonumber \\
&=&m\ns{0}\ot_{\mathcal{K}}t_2\ot_S\cdots \ot_S t_{n+1}\ot_S {m\ns{-1}}\triangleright( t_0 t_1)\nonumber \\
&=&\mathrm{t}_{n}(m\ot_{\mathcal{K}}t_0 t_1\ot_S t_2\ot_S\cdots \ot_S t_{n+1})\nonumber \\
&&\mathrm{t}_{n}s_{0}(m\ot_{\mathcal{K}}t_0\ot_S \cdots \ot_S t_{n+1}) \nonumber. \\
\end{eqnarray}
In the above, we use the left $\mathcal{K}$-coaction property of $M$ in the third equality and left $\mathcal{K}$-module algebra property of $T$ in the fifth equality. The relation $\mathrm{t}_n d_{0}= d_{n}$ is obvious. The other cocyclicity conditions can be easily verified.
\end{proof}

The cyclic cohomology of the preceding  cocyclic module is denoted by

$_\Kc{HC^\ast}(T,M)$ , and is called  Hopf cyclic cohomology of $T$ with coefficients in $M$ under the symmetry of $\mathcal{K}$.

To get a cyclic module we apply the duality procedure to  the previous cocyclic module as follows. Let $S$ be an algebra  over  $\mathbb{C}$, $\mathcal{K}$ a left $\times_{S}$-Hopf algebra, $T$ a left $\mathcal{K}$-module algebra and $M$ a right-left SAYD module over $\mathcal{K}$.  Let
\begin{equation}\label{cyclic-module-algebra}
_\Kc{C}_{n}(T, M)= M\otimes_{\mathcal{K}}T^{\otimes_{S}(n+1)}.
\end{equation}The following faces, degeneracies and cyclic map define a cyclic module structure on $_\Kc{C}_{\ast}(T, M)$.
\begin{align}\label{cyclic-module-algebra-operators}
\begin{split}
&\d_i (m\otimes_{\mathcal{K}} \td t)=m\otimes_{\mathcal{K}}t_0\otimes_S\cdots\ot_S t_i t_{i+1}\ot_{S}\cdots\ot_S t_n, \quad 0\leq i \leq n-1,\\
&\d_n (m\otimes_{\mathcal{K}} \td t)=m\ns{0}\ot_{\mathcal{K}} t_n (m\ns{-n}\triangleright t_0)\ot_S m\ns{-n+1}\triangleright t_1 \ot_S \cdots \ot_S m\ns{-1}\triangleright t_{n-1}, \\
&\s_i (m\otimes_{\mathcal{K}} \td t)=m\otimes_{\mathcal{K}}t_{0}\otimes_{S}\cdots\otimes_{S}t_{i}\otimes_{S}1\otimes_{S}t_{i+1}\otimes_S\cdots\otimes_{S}t_{n}, \quad 0\leq i \leq n, \\
&\tau_n(m\otimes_{\mathcal{K}} \td t)= m\ns{0} \ot_{\mathcal{K}} t_{n}\ot_{S} m\ns{-n}\triangleright t_0 \ot_{S} \cdots \ot_{S} m\ns{-1}\triangleright t_{n-1}.
\end{split}
\end{align}

 Now we recall the cyclic homology of a module coring with coefficients in a SAYD module under the symmetry of a $\times_{R}$-Hopf algebra.
  Let $\mathcal{B}$ be a right $\times_{R}$-Hopf algebra. A right $\mathcal{B}$-module coring  $C$  is an $R$-coring and right $\mathcal{B}$-module such that  $R$-bimodule structure of $C$ coincides with the one induced by $\Bc$. In addition,   it is assumed that the  counit $\varepsilon$ and comultiplication $\Delta$  are right $\mathcal{B}$ linear. That is  we consider the right $\mathcal{B}$-module structure of $R$ by $r\vartriangleleft b:= \varepsilon(\mathfrak{s}(r)b)$ and right $\mathcal{B}$ module structure of $C\otimes_{R}C$ is by the diagonal action. This means:
\begin{align}\label{mc1}
&\varepsilon_{C}(c\vartriangleleft b)= \varepsilon_{C}(c)\vartriangleleft b= \varepsilon_{\mathcal{B}}(\mathfrak{s}(\varepsilon_{C}(c))b),\\\label{mc2}
&\Delta_C(c\vartriangleleft b)= \Delta(c)\triangleleft b=c^{(1)} \vartriangleleft b^{(1)} \otimes_{R}c^{(2)} \vartriangleleft b^{(2)}.
\end{align}

We define
\begin{equation}
\widetilde{C}_{ \mathcal{B},n}(C, M)=C^{\otimes_{R}(n+1)}\otimes_{\mathcal{B}}M,
\end{equation}
and following operators on  $\widetilde{C}_{ \mathcal{B},n}(C, M)$.
\begin{align}
\begin{split}\label{module-ring-cyclic}
&\delta_{i}(\td c\otimes_{\mathcal{B}}m)= c_0\otimes_R\cdots\otimes_R \varepsilon_{C}(c_i)\otimes_R\cdots \ot_R c_n\ot_{\mathcal{B}}m, \quad 0\leq i \leq n,\\
& \sigma_i(\td c\otimes_{\mathcal{B}}m)= c_0\otimes_R\cdots\otimes_R \Delta_{C}(c_i)\otimes_R\cdots \ot_Rc_n\ot_{\mathcal{B}}m, \quad 0\leq i \leq n,\\
&\tau_{n}(\td c\otimes_{\mathcal{B}}m)= c_n \triangleleft m\ns{1}\otimes_R c_0\ot_R\cdots\ot_R c_{n-1}\ot_{\mathcal{B}}m\ns{0}.
\end{split}
\end{align}
Here $\td c=c_0\otimes_R\cdots\otimes_R c_n$ and  $C^{\ot_R (n+1)}$ is a right $\mathcal{B}$-module by diagonal action.
After some identifications the following theorem coincides with Proposition 2.19 in \cite{bs2}. For reader convenience we just check that the cyclic operator is  well defined.
\begin{proposition}
Let $\mathcal{B}$ be a right $\times_{R}$-Hopf algebra, $M$ a left-right SAYD module over $\mathcal{B}$ and $C$ a right $\mathcal{B}$-module  coring. Then the operators defined in {\rm \ref{module-ring-cyclic}} define a cyclic module structure on  $\widetilde{C}_{\mathcal{B},*}(C, M)$.
\end{proposition}
\begin{proof}
Since $C$ is a right $\mathcal{B}$-module coring we have, $(c\triangleleft b^{(1)})(r\triangleleft b^{(2)})=(cr)\triangleleft b,$ for $ c\in C$ and $ r\in R.$
Therefore $ (c_1\triangleleft b^{(1)})(\varepsilon(c_2)\triangleleft b^{(2)})=(c_1\varepsilon(c_2))\triangleleft b \quad \text{for} \quad c_1,c_2\in R$, shows that faces are well-defined. Also since comultiplication  is a right $\mathcal{B}$-module map, using \eqref{mc2} degeneracies are well-defined. The following calculation shows that the cyclic map is well-defined:
\begin{align*}
&\tau_{n}(c_0\otimes_R\cdots\ot_R c_n\otimes_{\mathcal{B}}b\triangleright m)= \\
&c_n \triangleleft (b\triangleright m)\ns{1}\ot_R c_0\cdots \ot_R c_{n-1}\ot_{\mathcal{B}}(b\triangleright m)\ns{0} =\\
&c_n \triangleleft b^{(2)}m\ns{1}{b^{(1)}}^{-}\ot_R c_0 \ot_R \cdots \ot_R c_{n-1}\ot_{\mathcal{B}} {b^{(1)}}^{+}\triangleright m\ns{0}=\\
&((c_n \triangleleft {b^{+}}^{(2)}m\ns{1}b^{-})\ot_R c_0\ot_R\cdots\ot_R c_{n-1})\triangleleft {b^{+}}^{(1)}\ot_{\mathcal{B}}m\ns{0} = \\
&((c_n \triangleleft {b^{+}}^{(2)}m\ns{1}b^{-}{b^{+}}^{(1)(1)})\ot_R c_0\triangleleft{b^{+}}^{(1)(2)}\ot_R\cdots\ot_R c_{n-1}\triangleleft{b^{+}}^{(1)(n+1)}\ot_{\mathcal{B}}m\ns{0}=  \\
&((c_n \triangleleft {b^{+}}^{(n+2)}m\ns{1}b^{-}{b^{+}}^{(1)})\ot_R c_0\triangleleft{b^{+}}^{(2)}\ot_R\cdots\ot_R c_{n-1}\triangleleft{b^{+}}^{(n+1)}\ot_{\mathcal{B}}m\ns{0}  =\\
&(c_n \triangleleft b^{(n+1)}m\ns{1})\ot_R c_0 \triangleleft b^{(1)}\ot_R\cdots\ot_R c_{n-1}\triangleleft b^{(n)}\ot_{\mathcal{B}}m\ns{0}  =\\
&\tau_{n}( c_0 \triangleleft b^{(1)}\ot_R\cdots\ot_R \ot_R c_n \triangleleft b^{(n+1)}\ot_{\mathcal{B}}m ) =\\
&\tau_{n}((c_0\otimes_R\cdots\otimes_R c_n)\triangleleft b\otimes_{\mathcal{B}}m). \\
\end{align*}
We use the AYD condition \eqref{bb} in the second equality, \cite[Lemma 2.14.(v)]{bs2} in the third equality, the diagonal action in the fourth equality and  \cite[Lemma 2.14.(i)]{bs2} in the sixth equality. We leave to the reader to check that $\d_i,\s_j$ and $\tau$ satisfy all conditions of a cyclic module.
\end{proof}

The cyclic homology of this cyclic module is  denoted by $\widetilde{HC}_{\Bc,*}(C, M)$.

In a special case, $\mathcal{B}$ is a right $\mathcal{B}$-module coring by multiplication map $\mathcal{B}\ot_R \mathcal{B}\longrightarrow \mathcal{B}.$
To synchronize with the Hopf cyclic complex of Hopf algebras we identify
\begin{equation}\label{cycli-module-X-Hopf}
\xymatrix{ \widetilde{C}_{\Bc,n}(\mathcal{B}, M)\ar[r]^{\varphi_n}&\mathcal{B}^{\ot_{R}n}\ot_{R^{op}} M,}
\end{equation}
where $\varphi$ defined by
\begin{align}\label{varphi}
&\varphi_n: \mathcal{B}^{\ot_R (n+1)}\ot_\Bc M\ra \mathcal{B}^{\ot_R n}\ot_{R^{op}} M,\\\notag
&\varphi_n(b_0\ot_R\cdots\ot_R  b_n\ot_{\mathcal{B}}m)=(b_0\ot_R\cdots\ot_R b_{n-1})\triangleleft b_{n}^{-}\ot_{R^{op}} b_{n}^{+}m,
\end{align}
\begin{proposition}
The map $\varphi_n$ defined in \eqref{varphi} is a well-defined isomorphism of vector spaces.
\end{proposition}
\begin{proof}
By \cite[Lemma 2.14. (iii)(ii)]{bs2}, the map $\varphi$ is well-defined. We define the inverse of map $\varphi$ by
$$\varphi_n^{-1}(b_1\ot_R\cdots\ot_R b_n\ot_{R^{op}} m)=b_1\ot_R\cdots\ot_R b_{n}\ot_R 1_{\mathcal{B}}\ot_{\mathcal{B}}m.$$ This map is well-defined  because for any $r^{o}\in R^{op}$ we have  $r^{o}m= \mathfrak{t}(r)\triangleright m$ and $(b_1\ot_R\cdots \ot_R b_n)\triangleleft r^{o}=b_1\mathfrak{t}(r)\ot_R\cdots \ot_R b_n.$ One can easily check that $\varphi \varphi^{-1}=Id$ by \cite[Lemma 2.14. (iv)]{bs2},  and $\varphi^{-1}\varphi=Id$ by \cite[Lemma 2.14.(i)]{bs2}.

\end{proof}
Therefore we transfer the cyclic structure of  $\widetilde{C}_{\Bc,*}(\Bc,M)$ to $\mathcal{B}^{\ot_R n}\ot_{R^{op}} M$. The resulting operators are recorded bellow.
\begin{align}\label{cyclic-module-X-Hopf-operators}
\begin{split}
&\delta_{i}(\td b\ot_{R^{op}} m)= b_1\otimes_R\cdots\otimes_R \varepsilon(b_{i+1})\otimes_R\cdots \ot_R b_n\ot_{R^{op}} m, \quad 0\leq i\leq n-1\\
&\delta_{n}(\td b\ot_{R^{op}} m)= (b_1\otimes_R \cdots \ot_Rb_{n-1})\triangleleft b_n^{-}\ot_{R^{op}} b_n^{+}m,\\
&\sigma_i(\td b \ot_{R^{op}} m)= b_1\otimes_R\cdots\otimes_R \Delta(b_{i+1})\otimes_R\cdots \ot_R b_n\ot_{R^{op}} m, \quad 0\leq i \leq n-1,\\
&\sigma_n(\td b\otimes_{R^{op}} m)= b_1\otimes_R\cdots\otimes_R b_n\ot_R 1_{\mathcal{B}} \ot_{R^{op}} m, \\
&\tau_n(\td b \ot_{R^{op}} m)= (m\ns{1}\ot_R b_1\ot_R\cdots\ot_R b_{n-1})\triangleleft b_{n}^{-}\ot_{R^{op}} b_n^{+}m\ns{0},
\end{split}
\end{align}
where $\td b= b_1\ot_R \cdots \ot_R b_n$.
%%%%%%%%%%%%%%%%%%%%%%%%%%%%%%%%%%%%%%%%%%%%%%%%%%%%%%%%%%%

\section{Equivariant Hopf Galois extensions}\label{S-3}
In this section we first define  equivariant Hopf Galois extensions. An equivariant Hopf Galois extension is a quadruple $(\Kc,\Bc,T,S)$ satisfying certain properties as stated in Definition \ref{galois}.   In the Subsection \ref{ss-3-1} we show that any equivariant Hopf Galois extension $(\Kc,\Bc,T,S)$ defines a functor from the category of SAYD modules over  $\Kc$ to the category of SAYD modules over $\Bc$ such that their Hopf cyclic complexes with  corresponding coefficients are isomorphic.
In the Subsection \ref{ss-3-2} we introduce an example of equivariant Hopf Galois extension and  explicitly illustrate the functor between the categories of SAYD modules.

Let  $\mathcal{B}$ be a right $\times_R$-Hopf algebra. We say the algebra $T$ is a  right $\mathcal{B}$-comodule algebra via $\Db:T\ra T\ot_R \Bc$ if  $T$ is $R$-bimodule and right $B$-comodule satisfying the following   conditions
\begin{enumerate}
\item[i)] $(t_1 r) t_2 = t_1( r t_2 ), $ (multiplication in $T$ is $R$-balanced).
\item[ii)] $1_T\ns{0}\ot_R 1_T\ns{1}= 1_T\ot_R 1_{\mathcal{B}},$
\item[iii)]$(tt')\ns{0} \ot_R (tt')\ns{1}= t\ns{0} t'\ns{0} \ot_R t\ns{1} t'\ns{1}.$
\end{enumerate}
The  coinvariants subalgebra $S\subseteq T$ is defined by
\begin{align}
S:=T^{\mathcal{B}}=\{ s\in T\mid \quad \text{ where} \quad \Db(s)= s\ot_R 1_{\mathcal{B}}\}.
\end{align}
\begin{definition}\label{extension} Let $R$  be an algebra over $\mathbb{C}$,  $\mathcal{B}$ a right $\times_R$-Hopf algebra, $T$ a right $\mathcal{B}$-comodule algebra, $S=T^{\mathcal{B}}$, $\mathcal{K}$ a left $\times_S$-Hopf algebra, and $T$ a left $\mathcal{K}$-module algebra. $T$ is called
 a $\mathcal{K}$-equivariant $\mathcal{B}$-Galois extension of $S$, if the canonical map
\begin{align}\label{galois}
can: T\ot_S T\longrightarrow T\ot_R \mathcal{B}, \quad t'\ot_S t\longmapsto t't\ns{0}\ot_R t\ns{1},
\end{align}
is bijective and the right coaction of $\mathcal{B}$  over $T$ is  $\mathcal{K}$-equivariant, i.e.,
\begin{align}
(kt)\ns{0}\ot_R (kt)\ns{1}= kt\ns{0}\ot_R t\ns{1}, \quad \quad k\in \mathcal{K} \quad t\in T.
\end{align}
\end{definition}
We denote a $\mathcal{K}$-equivariant $\mathcal{B}$-Galois extension described in the Definition \ref{extension} by $_{\mathcal{K}}T(S)^{\mathcal{B}}$.
One observes that by $\Kc:= B\ot B^{op}$ in the above definition we recover the ordinary Hopf Galois extension since  for any $\Bc$-comodule algebra $T$ the  Galois map \eqref{galois} is always $B\ot B^{op}$-equivariant.

%For example  since for $\mathcal{B}$-Galois extension $T(S)^{\mathcal{B}}$ over a right $\times_R$-Hopf algebra $\mathcal{B}$, described in \cite[Definition 2.21]{bs2} is a $\mathcal{K}$-equivariant $\mathcal{B}$-Galois extension, where $\mathcal{K}=B\ot B^{op}$ and $S=\mathcal{B}$.
\medskip

One defines  left $\mathcal{K}$-module structures on $T\ot_S T$ and $T\ot_R \mathcal{B}$, respectively, as follows.
\begin{align}
k\triangleright (t_1\ot_S t_2)= k^{(1)}t_1\ot_S k^{(2)}t_2, \quad \text{and} \quad k\triangleright(t\ot_R b)= kt\ot_R b.
\end{align}

The $\mathcal{K}$-equivariant condition of  $_{\mathcal{K}}T(S)^{\mathcal{B}}$ implies that
\begin{align}
can(k\triangleright (t'\ot_S t))= k\triangleright can(t'\ot_S t), \quad \quad k\in \mathcal{K} \quad t,t'\in T.
\end{align}

We denote  the inverse of the Galois map \eqref{galois} by the following index notation
\begin{align}
can^{-1}(1\ot_R b):= b\ns{-} \ot_S b\ns{+}.
\end{align}
One has the following properties for the maps $can$ and $can^{-1}$.
\begin{lemma}\label{inverse} Let $_{\mathcal{K}}T(S)^{\mathcal{B}}$ be a  $\mathcal{K}$-equivariant $\mathcal{B}$-Galois extension. Then the following properties hold.\\
\begin{enumerate}
\item[\rm i)]$(k\triangleright t)b\ns{-}\ot_S b\ns{+} = (k^{(1)}t)(k^{(2)}b\ns{-})\ot_S k^{(3)}b\ns{+}.$
\item[\rm ii)]$(b b')\ns{-}\ot_S (b b')\ns{+}= b'\ns{-}b\ns{-}\ot_S b\ns{+}b'\ns{+}.$
\item[\rm iii)] $b\ns{-}\ns{0}\ot_R b\ns{-}\ns{1} \ot_S b\ns{+}= b^{+}\ns{-}\ot_R b^{-}\ot_S b^{+}\ns{+}.$
\item[\rm iv)] $b\ns{-}\ot_S b\ns{+}\ns{0}\ot_R b\ns{+}\ns{1}=b^{(1)}\ns{-}\ot_S b^{(1)}\ns{+}\ot_R b^{(2)}.$
\item[\rm v)] $b\ns{-}b\ns{+}\ns{0}\ot_R b\ns{+}\ns{1}= 1\ot_R b.$
\item[\rm vi)] $t\ns{0}t\ns{1}\ns{-}\ot_S t\ns{1}\ns{+}=1\ot_S t.$
\item[\rm vii)] $1_{\mathcal{B}}\ns{-}\ot_S 1_{\mathcal{B}}\ns{+}= 1_{T}\ot_S 1_T.$
\item[\rm viii)] $\mathfrak{t}(r)\ns{-}\ot_S \mathfrak{t}(r)\ns{+}= r\triangleright 1_T \ot_S 1_T.$
\item[\rm ix)] $\mathfrak{s}(r)\ns{-}\ot_S \mathfrak{s}(r)\ns{+}= 1_T \ot_S 1_T\triangleleft r.$
\end{enumerate}
\end{lemma}
\begin{proof}
One easily sees that  (i)  is equivalent to $\mathcal{K}$-equivariant  property of the map $can^{-1}$. One can find  a proof of the rest in \cite[pages 268-270]{bs2}.

\end{proof}
%%%%%%%%%%%%%%%%%%%%%%%%%%%%%%%%%%%%%%%%%%%%%%%%%%%%%%%%%%%%%%%%%%%%%%%%%%%%%%%%%%%%%%%%%%%%%%%%%%%%%%%%%%
\subsection{Equivariant Hopf Galois extension as a functor}{\label{ss-3-1}

Let $_\Kc T(S)^\Bc$ be a $\Kc$-equivariant $\Bc$-Galois extension, and $M$ be a right-left SAYD module over $\Kc$.  We let  $\mathcal{B}$ act on $\widetilde{M}:= M\ot_\Kc T$ on the left  by
\begin{equation}\label{action}
b\triangleright (m\ot_{\mathcal{K}} t)= m\ns{0}\ot_{\mathcal{K}} b\ns{+}(m\ns{-1}\triangleright (t b\ns{-})).
\end{equation}
 We also let $\Bc$ coact on $\widetilde{M}$ from the right by
\begin{equation}\label{coaction}
(m\ot_{\mathcal{K}} t)\longmapsto m\ot_{\mathcal{K}} t\ns{0} \ot_R t\ns{1}.
\end{equation}
\begin{theorem}\label{main} Let $R$  be an algebra over $\mathbb{C}$,  $\mathcal{B}$ a right $\times_R$-Hopf algebra, $T$ a right $\mathcal{B}$-comodule algebra, $S=T^{\mathcal{B}}$, $\mathcal{K}$ a left $\times_S$-Hopf algebra, $T$ a left $\mathcal{K}$-module algebra and $M$ be a right-left SAYD module  over $\mathcal{K}$. If $_{\mathcal{K}}T(S)^{\mathcal{B}}$ is a $\mathcal{K}$-equivariant $\mathcal{B}$-Galois extension, then $\widetilde{M}:= M\ot_{\mathcal{K}} T$ is a left-right SAYD module over $\mathcal{B}$ by the action and coaction defined in \eqref{action} and \eqref{coaction}.
\end{theorem}
\begin{proof}

The above coaction is obviously well-defined  by the $\mathcal{K}$-equivariant property of the coaction of $\mathcal{B}$ over $T$. We shall  show that the action \eqref{action} is well-defined. We show that
\begin{equation*}
b\triangleright (m\triangleleft k\ot_{\mathcal{K}} t)= b\triangleright(m\ot_{\mathcal{K}} k\triangleright t), \quad  \quad b\in\mathcal{B}, \quad m\in M, \quad k\in \mathcal{K}, \quad t\in T.
\end{equation*}
Indeed,
\begin{align*}
b\triangleright (m\triangleleft k\ot_{\mathcal{K}} t)&= (mk)\ns{0}\ot_{\mathcal{K}} b\ns{+}((mk)\ns{-1}\triangleright (tb\ns{-})) \\
&=m\ns{0}{k^{(2)}}^{-}\ot_{\mathcal{K}} b\ns{+} ({k^{(2)}}^{+}m\ns{-1}k^{(1)}\triangleright (tb\ns{-}))  \\
&=m \ns{0}\ot_{\mathcal{K}} {k^{(2)}}^{-}\triangleright(b\ns{+} ({k^{(2)}}^{+}m\ns{-1}k^{(1)}\triangleright (tb\ns{-})))  \\
&=m \ns{0}\ot_{\mathcal{K}} ({{k^{(2)}}^{-}}^{(1)}b\ns{+}) ({{k^{(2)}}^{-}}^{(2)}{k^{(2)}}^{+}m\ns{-1}k^{(1)}\triangleright (tb\ns{-})))  \\
&=m\ns{0} \ot_{\mathcal{K}} (k^{(2)}\triangleright b\ns{+}) (m\ns{-1}k^{(1)}\triangleright (tb\ns{-}))  \\
&=m\ns{0} \ot_{\mathcal{K}} (k^{(3)}\triangleright b\ns{+}) (m\ns{-1}\triangleright((k^{(1)}\triangleright t) (k^{(2)}\triangleright b\ns{-})))  \\
&=m\ns{0} \ot_{\mathcal{K}}b\ns{+} (m\ns{-1}\triangleright((k\triangleright t) b\ns{-}))  \\
&=b\triangleright(m\ot_{\mathcal{K}} k\triangleright t).
\end{align*}
We use the SAYD condition \eqref{b} for $M$ over $\mathcal{K}$ in the second equality, the $\mathcal{K}$-module algebra property of $T$ in the fourth equality,  the Lemma \ref{lem}(i) in the fifth equality, the $\mathcal{K}$-module algebra property of $T$ in the sixth equality and  the Lemma \ref{inverse}(i)  in the seventh equality. Next we show that the action \eqref{action} is $S$-balanced both in $b\ns{-}\ot_S b\ns{+}$ and $m\ns{-1}\ot_S m\ns{0}$. That is for any $q,p\in S, m\in M, k\in \mathcal{K}$ and $t_1, t_2, t_3 \in T$, we should show
\begin{equation}
 pm\ot_{\mathcal{K}} t_1(k\triangleright t_2 t_3)=m\ot_{\mathcal{K}}t_{1}(\mathfrak{t}(p)k\triangleright t_2 t_3),
\end{equation}
which is obvious because
\begin{eqnarray}
  && pm\ot_{\mathcal{K}} t_1(k\triangleright t_2 t_3)= m\triangleleft \mathfrak{t}(p)\ot_{\mathcal{K}} t_1(k\triangleright t_2 t_3)\nonumber \\
  &=& m\ot_{\mathcal{K}} \mathfrak{t}(p)\triangleright t_1(k\triangleright t_2 t_3)= m\ot_{\mathcal{K}}  t_1(\mathfrak{t}(p)k\triangleright t_2 t_3), \nonumber
\end{eqnarray}
where $\mathfrak{t}$ is the target map of the left bialgebroid $\mathcal{K}$. In the above we use $pm=  m\triangleleft \mathfrak{t}(p)$ and $\Delta(\mathfrak{t}(p))=1_{\mathcal{K}}\ot_S \mathfrak{t}(p)$.
This allows then to substitute $k\ot_S m= n\ns{-1}\ot_S n\ns{0}$ for any $n\in M$. To finish the argument that the action \eqref{action} is $S$-balanced we should also show
\begin{equation}
  m\ns{0}\ot_{\mathcal{K}} t_1 (m\ns{-}\triangleright t_2 t_3 q)= m\ns{0}\ot_{\mathcal{K}} q t_1 (m\ns{-}\triangleright t_2 t_3 ),
\end{equation}
which is clear because
\begin{eqnarray}
  && m\ns{0}\ot_{\mathcal{K}} t_1 (m\ns{-1}\triangleright t_2 t_3 q)= m\ns{0}\ot_{\mathcal{K}} t_1 (m\ns{-1}\mathfrak{t}(q)\triangleright t_2 t_3 )\nonumber \\
  &=& m\ns{0}q\ot_{\mathcal{K}} t_1 (m\ns{-1}\triangleright t_2 t_3 )=  m\ns{0}\triangleleft \mathfrak{s}(q)\ot_{\mathcal{K}} t_1 (m\ns{-1}\triangleright t_2 t_3 )\nonumber\\
  &=&m\ns{0}\ot_{\mathcal{K}} \mathfrak{s}(q)\triangleright t_1 (m\ns{-1}\triangleright t_2 t_3) =m\ns{0}\ot_{\mathcal{K}} q t_1 (m\ns{-1}\triangleright t_2 t_3), \nonumber
\end{eqnarray}
where $\mathfrak{s}$ is the source map of the left bialgebroid $\mathcal{K}$. In the above we use \eqref{abc} in the second equality. Thus it makes sense to substitute $t_3\ot_S t_1= b\ns{-}\ot_S b\ns{+}.$
Next we show that \eqref{action} is associative. For this, let $g,h\in \mathcal{B}$ and $m\ot_{\mathcal{K}} t \in \widetilde{M}$.
\begin{align*}
&g\triangleright(h\triangleright(m\ot_{\mathcal{K}} t))=g\triangleright(m\ns{0}\ot_{\mathcal{K}} h\ns{+}(m\ns{-1}\triangleright ( th\ns{-})))\\
&=m\ns{0}\ns{0}\ot_{\mathcal{K}} g\ns{+}\{m\ns{0}\ns{-1}\triangleright[h\ns{+}(m\ns{-1}\triangleright(th\ns{-}))g\ns{-}]\}\\
&=m\ns{0}\ot_{\mathcal{K}} g\ns{+}\{{m\ns{-1}}^{(2)}\triangleright[h\ns{+}({m\ns{-1}}^{(1)}\triangleright(th\ns{-}))g\ns{-}]\}\\
&=m\ns{0}\ot_{\mathcal{K}} g\ns{+}\{{m\ns{-1}}^{(2)(1)}\triangleright [h\ns{+}({m\ns{-1}}^{(1)}\triangleright(th\ns{-}))]\} ({m\ns{-1}}^{(2)(2)} \triangleright g\ns{-})\\
&=m\ns{0}\ot_{\mathcal{K}} g\ns{+}\{{m\ns{-1}}^{(2)(1)}\triangleright [h\ns{+}({m\ns{-1}}^{(1)(1)}\triangleright t)({m\ns{-1}}^{(1)(2)}\triangleright h\ns{-})]\} \\ &({m\ns{-1}}^{(2)(2)} \triangleright g\ns{-})\\
&=m\ns{0}\ot_{\mathcal{K}} g\ns{+}\{{m\ns{-1}}^{(3)}\triangleright [ h\ns{+}({m\ns{-1}}^{(1)}\triangleright t)({m\ns{-1}}^{(2)}\triangleright h\ns{-})]\}\\
&({m\ns{-1}}^{(4)} \triangleright g\ns{-})\\
%&=m\ns{0}\ot_{\mathcal{K}} g\ns{+}\{{m\ns{-1}}^{(2)(1)}\triangleright [h\ns{+}({m\ns{-1}}^{(1)}\triangleright t)({m\ns{-1}}^{(2)}\triangleright h\ns{-})]\}({m\ns{-1}}^{(2)(2)} \triangleright g\ns{-})\\
%&=m\ns{0}\ot_{\mathcal{K}} g\ns{+}\{{m\ns{-1}}^{(3)}\triangleright [h\ns{+}({m\ns{-1}}^{(1)}\triangleright t)({m\ns{-1}}^{(2)}\triangleright h\ns{-})]\}({m\ns{-1}}^{(4)} \triangleright g\ns{-})\\
&=m\ns{0}\ot_{\mathcal{K}} g\ns{+}({m\ns{-1}}^{(3)}\triangleright  h\ns{+})[({m\ns{-1}}^{(4)}\triangleright({m\ns{-1}}^{(1)}\triangleright t)({m\ns{-1}}^{(2)}\triangleright h\ns{-})]\\
&({m\ns{-1}}^{(5)} \triangleright g\ns{-})\\
&=m\ns{0}\ot_{\mathcal{K}} g\ns{+}h\ns{+}[({m\ns{-1}}^{(2)}\triangleright({m\ns{-1}}^{(1)}\triangleright t)h\ns{-}]({m\ns{-1}}^{(3)}\triangleright g\ns{-})\\
&=m\ns{0}\ot_{\mathcal{K}} g\ns{+}h\ns{+}[{m\ns{-1}}^{(2)}\triangleright({m\ns{-1}}^{(1)}\triangleright t)h\ns{-}g\ns{-})]\\
&=m\ns{0}\ns{0}\ot_{\mathcal{K}}(gh)\ns{+}[m\ns{0}\ns{-1}\triangleright(m\ns{-1}\triangleright t)(gh)\ns{-}]\\
&=(gh)\triangleright(m\ns{0}\ot_{\mathcal{K}} m\ns{-1}\triangleright t)\\
&=(gh)\triangleright(m\ot_{\mathcal{K}} t).
\end{align*}
We apply  the coaction property in the second equality, the $\mathcal{K}$-module algebra property of $T$ in the third, fourth, sixth and eighth equalities, the $\mathcal{K}$-equivariance property of Lemma \eqref{inverse}(i) in the seventh equality  and the stability condition in the last equality. Also unitality of the action \eqref{action} is obvious by Lemma \eqref{inverse}(vii) and stability condition.
Now we verify the AYD condition for $\widetilde{M}$ over $\mathcal{B}$. Let $b\in B$, $t\in T$ and $m\in M$.
\begin{align*}
&(b\triangleright(m\ot_{\mathcal{K}} t))\ns{0}\ot_R (b\triangleright(m\ot_{\mathcal{K}} t))\ns{1}\\
&=(m\ns{0}\ot_{\mathcal{K}} b\ns{+}(m\ns{-1}\triangleright (t b\ns{-})))\ns{0} \ot_R \\
&(m\ns{0}\ot_{\mathcal{K}} b\ns{+}(m\ns{-1}\triangleright (t b\ns{-})))\ns{1}\\
&=m\ns{0}b\ns{+}\ns{0}(m\ns{-1}\triangleright (t\ns{0} b\ns{-}\ns{0})) \ot_R b\ns{+}\ns{1} t\ns{1} b\ns{-}\ns{1}\\
&=m\ns{0}{b^{(1)}}\ns{+}(m\ns{-1}\triangleright (t\ns{0} {b^{(1)}}\ns{-}\ns{0})) \ot_R b^{(2)} t\ns{1} {b^{(1)}}\ns{-}\ns{1}\\
&=m\ns{0}\ot_{\mathcal{K}} {b^{(1)}}^{+}\ns{+}(m\ns{-1}\triangleright (t\ns{0}{b^{(1)}}^{+}\ns{-}))\ot_R b^{(2)} t\ns{1}{b^{(1)}}^{-}\\
&={b^{(1)}}^{+}\triangleright(m\ot_{\mathcal{K}}t)\ns{0}\ot_R b^{(2)}(m\ot_{\mathcal{K}} t)\ns{1}{b^{(1)}}^{-}.\\
\end{align*}

We use the $\mathcal{K}$-equivariant property of the coaction of $\mathcal{B}$ on $T$  and the $\mathcal{B}$-comodule algebra property of  $T$ in the second equality, the Lemma \ref{inverse}(iv) in the third equality and  the Lemma \ref{inverse}(iii) in the fourth equality.

 The stability condition of $\widetilde{M}$ over $\mathcal{B}$ follows from
\begin{align*}
&(m\ot_{\mathcal{K}} t)\ns{1}(m\ot_{\mathcal{K}} t)\ns{0}= t\ns{1}(m\ot_{\mathcal{K}} t\ns{0})\\
&= m\ns{0}\ot_{\mathcal{K}} (t\ns{1}\ns{+}(m\ns{-1}\triangleright (t\ns{0}t\ns{1}\ns{-})))\\
&=m\ns{0}\ot_{\mathcal{K}} t(m\ns{-1}\triangleright 1_T)=m\ns{0}\ot_{\mathcal{K}} t(\mathfrak{t}(\varepsilon(m\ns{-1}))\triangleright 1_T)\\
&=m\ns{0}\ot_{\mathcal{K}} \mathfrak{t}(\varepsilon(m\ns{-1}))\triangleright t= m\ns{0}\triangleleft \mathfrak{t}(\varepsilon(m\ns{-1}))\ot_{\mathcal{K}} t\\
&=m\ot_{\mathcal{K}} t,
\end{align*}
Lemma \ref{inverse}~ (vi) is  applied  in the third equality, the module algebra property in fourth equality and the counitality  of the coaction on $M$ is used in the last equality.
\end{proof}
Let us denote the category of SAYD modules over a left $\times$-Hopf algebra $\Kc$ by $^\Kc SAYD_\Kc$. Its objects   are all right-left  SAYD modules over $\Kc$ and its morphisms  are all $\Kc$- linear-colinear   maps.  Similarly one denotes by  $_\Bc SAYD^\Bc$  the category of left-right SAYD modules of a right  $\times$-Hopf algebra $\Bc$. We see that Theorem \ref{main} amounts  to the object map of a functor $\digamma$ from $^\Kc SAYD_\Kc$ to $_\Bc SAYD^\Bc$. Let $\phi: M\ra N$ be a  morphism in $^\Kc SAYD_\Kc$ we define $\digamma(\phi)$ to be $\phi\ot_\Kc \Id_{T}$.
\begin{proposition}\label{diagmma}
The assignment $\digamma:\;\;^\Kc SAYD_\Kc\ra\;\; _\Bc SAYD^\Bc$ defines a covariant functor.
\end{proposition}
\begin{proof}
 Using Theorem \ref{main} we see that $\digamma$ is an object map. We now prove that it is also a morphism map.  Let $\phi: M\ra N$ be a  morphism in $^\Kc SAYD_\Kc$.  The colinearity of $\digamma(\phi)$ is obvious. The following shows that $\digamma(\phi)$ is also linear,
\begin{align*}
& (\phi\ot Id)(b\triangleright(m\ot_{\mathcal{K}}t))=(\phi\ot Id)(m\ns{0}\ot_{\mathcal{K}} b\ns{+}(m\ns{-1}\triangleright (t b\ns{-})))=\\
&\phi(m\ns{0})\ot_{\mathcal{K}} b\ns{+}(m\ns{-1}\triangleright (t b\ns{-}))= \phi(m)\ns{0}\ot_{\mathcal{K}} b\ns{+}(\phi(m)\ns{-1}\triangleright (t b\ns{-}))\\
&= b\triangleright (\phi(m)\ot t)=b\triangleright((\phi\ot Id)(m\ot_{\mathcal{K}}t)) .
\end{align*}
We use the comodule map property of $\phi$ in the third equality. The covariant property of the functor is obvious.

\end{proof}
Let us recall that $_\Kc{C}_n(T,M)= M\ot_\Kc T^{\ot_S n+1}$ and   $\nobreak{C_{\Bc,n}(B,\widetilde{M}):= B^{\ot_R n}\ot_{R^{op}} \widetilde{M} }$ are the cyclic modules defined in \eqref{cyclic-module-algebra-operators} and \eqref{cyclic-module-X-Hopf-operators} respectively. Now we define the following maps

\begin{align}
\omega_n: M\ot_{\mathcal{K}} T^{\ot_S (n+1)}\longrightarrow \mathcal{B}^{\ot_R n}\ot_{R^{op}} \widetilde{M},
\end{align}
given by
\begin{align*}
\omega_n(m\ot_{\mathcal{K}} t_0\ot_S\cdots \ot_S t_n)&=t_1\ns{1}t_2\ns{1}\cdots t_n\ns{1}\ot_R t_2\ns{2}\cdots t_n\ns{2}\ot_R\\
 &\cdots \ot_R t_n\ns{n}\ot_{R^{op}} (m\ot_{\mathcal{K}} t_0 t_1\ns{0}\cdots t_n\ns{0}),
\end{align*}
and
\begin{align*}
&\omega_n^{-1}(b_1\ot_R\cdots\ot_R b_{n+1}\ot_{R^{op}}m\ot_K t)= m\ot_K  tb_1\ns{-}\ot_S b_1\ns{+}b_2\ns{-}\ot_S\\
    & b_2\ns{+}b_3 \ns{-}\ot_S \cdots \ot_S b_{n-1}\ns{+}b_{n}\ns{-}\ot_S b_n\ns{+}.
\end{align*}
\begin{theorem}\label{general}
Let $_{\mathcal{K}}T(S)^{\mathcal{B}}$ be  a $\mathcal{K}$-equivariant $\mathcal{B}$-Galois extension, and  $M$ be a right-left SAYD module over $\Kc$.  Then  $\om_\ast$ defines an isomorphism of cyclic modules between   $_\Kc{C}_{\ast}(T,M)$ and $\widetilde{C}_{\mathcal{B},\ast}(\mathcal{B},\widetilde{M})$,  which are defined in \eqref{cyclic-module-algebra-operators} and \eqref{cyclic-module-X-Hopf-operators} respectively. Here   $\widetilde M:= M\ot_\Kc T$ is the left-right SAYD module over $\Bc$ introduced in Theorem \ref{main}.
\end{theorem}
\begin{proof} The map $\omega_n$ is the composition of the following maps

   \begin{equation}\begin{CD}
     M\ot_K T^{\ot_S(n+1)}@> \alpha_n >> \displaystyle \widetilde{M}\ot_R \mathcal{B}^{\ot_R n}
    @>\psi_n >>    \mathcal{B}^{\ot_R n}\ot_{R^{op}}\widetilde{M}\end{CD}.
    \end{equation}
    Here $\a_\ast$
 and $\psi_\ast$ are given by
        \begin{align*}
& \alpha_n(m\ot_{\mathcal{K}} t_0\ot_S\cdots \ot_S t_n)= m\ot_{\mathcal{K}} t_0 t_1\ns{0}\cdots t_n\ns{0}\ot_R\\
 &t_1\ns{1}t_2\ns{1}\cdots t_n \ns{1} \ot_R  \cdots\ot_R t_{n-1}\ns{n-1}t_n\ns{n-1}\ot_R t_n\ns{n},\\
  & \alpha_n^{-1}(m\ot_{\mathcal{K}} t\ot_R b_1\ot_R\cdots \ot_R b_n)= m\ot_{\mathcal{K}}  tb_1\ns{-}\ot_S\\
   &b_1\ns{+}b_2\ns{-}\ot_S b_2\ns{+}b_3 \ns{-}\ot_S \cdots \ot_S b_{n-1}\ns{+}b_{n}\ns{-}\ot_S b_n\ns{+}.\\
    &\\
    &\psi_n(m\ot_{\mathcal{K}} t \ot_R b_1\ot_R\cdots \ot_R b_n)= b_1\ot_R\cdots \ot_R b_n \ot_{R^{op}} m\ot_{\mathcal{K}} t, \\
    &\psi_n^{-1}(b_1\ot_R\cdots \ot_R b_n\ot_{R^{op}} m\ot_{\mathcal{K}} t)= m\ot_{\mathcal{K}} t \ot_R b_1\ot_R\cdots \ot_R b_n.
    \end{align*}

   One can easily use the $\mathcal{K}$-equivariant property of the coaction of $\mathcal{B}$ over $T$ and the $\mathcal{K}$-module algebra property of $T$ to see that the map $ \alpha_n$ is well-defined. Also using $\Delta(\mathfrak{s}(r))=1_{\mathcal{B}}\ot_{R} \mathfrak{s}(r)$ and $r\triangleright (b_1\ot_R\cdots\ot_R b_n)= b_1\mathfrak{s}(r)\ot_R\cdots\ot_R b_n$, the morphism $\psi_n$ is obviously  well-defined. We show that $\omega$ is a map of cyclic modules.
  By mutiplicity of the coproduct of $\mathcal{B}$, right $\mathcal{B}$-comodule algebra property of $T$  and the right $\times_R$-Hopf algebra property $b^{(2)}\mathfrak{t}(\varepsilon(\mathfrak{t}(r)b^{(1)}))= \mathfrak{t}(r)b$, the maps $\d_i, 0\leq i \leq n-1$ commute with $\omega$ as follows.

 \begin{align*}
& \d_i (\omega_n(m\ot_{\mathcal{K}} t_0\ot_S\cdots \ot_S t_n))= \d_i(t_1\ns{1}t_2\ns{1}\cdots t_n\ns{1}\ot_R t_2\ns{2}\\
&\cdots t_n\ns{2}\ot_R \cdots \ot_R t_n\ns{n}\ot_{R^{op}} (m\ot_{\mathcal{K}} t_0 t_1\ns{0}\cdots t_n\ns{0}))\\
&=t_1\ns{1}t_2\ns{1}\cdots t_n\ns{1}\ot_R \cdots \varepsilon(t_{i+1}\ns{i+1}\cdots t_n\ns{i+1})\ot_R\\
 &\cdots \ot_R t_n\ns{n}\ot_{R^{op}} (m\ot_{\mathcal{K}} t_0 t_1\ns{0}\cdots t_n\ns{0})\\
&=t_1\ns{1}t_2\ns{1}\cdots t_n\ns{1}\ot_R \cdots \ot_R t_{i}\ns{i}\cdots t_n\ns{i}\ot_R t_{i+2}\ns{i+1}\cdots t_n\ns{i+1}\ot_R\\
&\cdots \ot_R t_n\ns{n-1}\ot_{R^{op}} (m\ot_{\mathcal{K}} t_0 t_1\ns{0}\cdots t_n\ns{0})\\
&=t_1\ns{1}t_2\ns{1}\cdots t_n\ns{1}\ot_R\cdots \ot_R (t_{i}t_{i+1})\ns{i}t_{i+2}\ns{i}\cdots t_n\ns{i}\ot_R \\
&t_{i+2}\ns{i+1}\cdots t_n\ns{i+1}\ot_R\cdots \ot_R t_n\ns{n-1}\ot_{R^{op}}\\
 &(m\ot_{\mathcal{K}} t_0 t_1\ns{0}\cdots (t_i t_{i+1})\ns{0} t_{i+2}\ns{0}\cdots t_n\ns{0})\\
&=\omega_{n-1}(m\ot_{\mathcal{K}} t_0\ot_S\cdots \ot_S t_i t_{i+1}\ot_S \cdots \ot_S t_n)\\
&=\omega_{n-1}(\d_i(m\ot_{\mathcal{K}} t_0\ot_S\cdots \ot_S t_n))).
\end{align*}

Next we show that $\omega$ commutes with the last face morphism.

 \begin{align*}
& \d_n (\omega_n(m\ot_{\mathcal{K}} t_0\ot_S\cdots \ot_S t_n))=\d_n(t_1\ns{1}t_2\ns{1}\cdots t_n\ns{1}\ot_R\\
&t_2\ns{2}\cdots t_n\ns{2}\ot_R \cdots \ot_R t_n\ns{n}\ot_{R^{op}} (m\ot_{\mathcal{K}} t_0 t_1\ns{0}\cdots t_n\ns{0}))\\
&=(t_1\ns{1}\cdots t_n\ns{1}\ot_R \cdots \ot_Rt_{n-1}\ns{n-1}t_n\ns{n-1})\triangleleft {t_{n}\ns{n}}^{-}\ot_{R^{op}} \\
&{t_{n}\ns{n}}^{+}\triangleright(m\ot_{\mathcal{K}} t_0t_1\ns{0}\cdots t_n\ns{0})\\
&=(t_1\ns{1}\cdots t_{n-1}\ns{1}\ot_R  \cdots \ot_Rt_{n-1}\ns{n-1})\triangleleft t_n\ns{1}{t_{n}\ns{2}}^{-}\ot_{R^{op}} \\
&{t_{n}\ns{2}}^{+}\triangleright(m\ot_{\mathcal{K}} t_0t_1\ns{0}\cdots t_n\ns{0})\\
&=t_1\ns{1}\cdots t_{n-1}\ns{1}\ot_R \cdots \ot_Rt_{n-1}\ns{n-1}\ot_{R^{op}}\\
& t_n\ns{1}\triangleright(m\ot_{\mathcal{K}} t_0t_1\ns{0}\cdots t_n\ns{0})\\
&=t_1\ns{1}\cdots t_{n-1}\ns{1}\ot_R \cdots \ot_Rt_{n-1}\ns{n-1}\ot_{R^{op}} \\
&m\ns{0}\ot_{\mathcal{K}} {t_n\ns{1}}\ns{+}(m\ns{-1}\triangleright (t_0t_1\ns{0}\cdots t_n\ns{0}{t_n\ns{1}}\ns{-}))\\
&=t_1\ns{1}\cdots t_{n-1}\ns{1}\ot_R \cdots \ot_Rt_{n-1}\ns{n-1}\ot_{R^{op}} \\
&m\ns{0}\ot_{\mathcal{K}} t_n(m\ns{-1}\triangleright (t_0t_1\ns{0}\cdots t_{n-1}\ns{0}))\\
&=t_1\ns{1}\cdots t_{n-1}\ns{1}\ot_R \cdots \ot_Rt_{n-1}\ns{n-1}\ot_{R^{op}} \\
&m\ns{0}\ot_{\mathcal{K}} t_n(m\ns{-n}\triangleright t_0)(m\ns{-n+1}\rt t_1\ns{0})\cdots (m\ns{-1}\rt t_{n-1}\ns{0})\\
&=\omega_{n-1}(m\ns{0}\ot_{\mathcal{K}} t_n (m\ns{-n}\rt t_0)\ot_S m\ns{-n+1}\rt t_1\ot_S \cdots \ot_S m\ns{-1}\rt t_{n-1})\\
&=\omega_{n-1}\d_{n}(m\ot_{\mathcal{K}} t_0\ot_S\cdots \ot_S t_n).
  \end{align*}
 We apply the right diagonal action of $\mathcal{B}$ over $\mathcal{B}^{\ot_R (n-1)}$ in the third equality, the lemma \cite[2.14 .(ii)]{bs2} in the fourth equality, the action \eqref{action} in the fifth equality, the Lemma \eqref{inverse}(vi) in the sixth equality, the $\mathcal{K}$-module algebra property of $T$ in the seventh equality and $\mathcal{K}$-equivariant property of the coaction of $\mathcal{B}$ over $T$ in the eighth equality. By multiplicity of coproduct of  $\mathcal{B}$ the commutativity of $\omega$ and $\sigma_i, \quad 0\leq i \leq n-1,$ can be verified as follows.
   \begin{align*}
 &\sigma_{i}(\omega_n(m\ot_{\mathcal{K}} t_0\ot_S\cdots \ot_S t_n))  \\
&=\sigma_i(t_1\ns{1}t_2\ns{1}\cdots t_n\ns{1}\ot_R t_2\ns{2}\cdots t_n\ns{2}\ot_R\\
 &\cdots \ot_R t_n\ns{n}\ot_{R^{op}} (m\ot_{\mathcal{K}} t_0 t_1\ns{0}\cdots t_n\ns{0}))\\
 &=t_1\ns{1}t_2\ns{1}\cdots t_n\ns{1}\ot_R (t_i\ns{i+1}\cdots t_n\ns{i+1})^{(1)}\ot_R (t_{i+1}\ns{i+1}\cdots t_n\ns{i+1})^{(2)}\ot_R \cdots \ot_R\\
 &t_n\ns{n}\ot_{R^{op}} (m\ot_{\mathcal{K}} t_0 t_1\ns{0}\cdots t_n\ns{0})\\
&=t_1\ns{1}t_2\ns{1}\cdots t_n\ns{1}\ot_R t_{i+1}\ns{i+1}\cdots t_n\ns{i+1}\ot_R t_{i+1}\ns{i+2}\cdots t_n\ns{i+2} \cdots \ot_R\\
 &t_n\ns{n+1}\ot_{R^{op}} (m\ot_{\mathcal{K}} t_0 t_1\ns{0}\cdots t_n\ns{0})\\
&=\omega_{n+1}(m\ot_{\mathcal{K}} t_0\ot_S\cdots t_i\ot_S 1_T\ot_S t_{i+1}\ot_S\cdots \ot_S t_n \ot_S )\\
&=\omega_{n+1}( \sigma_{i}(m\ot_{\mathcal{K}} t_0\ot_S\cdots \ot_S t_n)).
   \end{align*}

Using Lemma \cite[2.14.(iv)]{bs2}, the commutativity of $\omega$ with the last degeneracy follows from

   \begin{align*}
& \sigma_{n}(\omega_n(m\ot_{\mathcal{K}} t_0\ot_S\cdots \ot_S t_n))~~~~~~~~~~~~~~~~~~~~~~~~~~~~~~~~~~~~~~~~~~~~~~~~~~~~~~~~~~~~ \\
&=\sigma_n(t_1\ns{1}t_2\ns{1}\cdots t_n\ns{1}\ot_R t_2\ns{2}\cdots t_n\ns{2}\ot_R\\
&\cdots \ot_R t_n\ns{n}\ot_{R^{op}} (m\ot_{\mathcal{K}} t_0 t_1\ns{0}\cdots t_n\ns{0}))\\
&=t_1\ns{1}t_2\ns{1}\cdots t_n\ns{1}\ot_R t_2\ns{2}\cdots t_n\ns{2}\ot_R \cdots \ot_R\\
&t_n\ns{n}\ot_R 1_T \ot_{R^{op}} (m\ot_{\mathcal{K}} t_0 t_1\ns{0}\cdots t_n\ns{0})\\
&=\omega_{n+1}(m\ot_{\mathcal{K}} t_0\ot_S\cdots\ot_S t_n \ot_S 1_T)\\
&=\omega_{n+1}( \sigma_{n}(m\ot_{\mathcal{K}} t_0\ot_S\cdots \ot_S t_n)).
   \end{align*}
 Finally, we show that $\omega$  commutes with the cyclic maps.
  \begin{align*}
& \tau_n (\omega_n(m\ot_{\mathcal{K}} t_0\ot_S\cdots \ot_S t_n))\\
&=\tau_n(t_1\ns{1}t_2\ns{1}\cdots t_n\ns{1}\ot_R t_2\ns{2}\cdots t_n\ns{2}\ot_R\\
 &\cdots \ot_R t_n\ns{n}\ot_{R^{op}} (m\ot_{\mathcal{K}} t_0 t_1\ns{0}\cdots t_n\ns{0}))\\
&=(t_0\ns{1}t_1\ns{0}\ns{1}\cdots t_n\ns{0}\ns{1}\ot_R t_1\ns{1}\cdots t_n\ns{1}\ot_R\\
 &\cdots \ot_Rt_{n-1}\ns{n-1}t_n\ns{n-1})\triangleleft {t_{n}\ns{n}}^{-}\ot_{R^{op}} \\
&{t_{n}\ns{n}}^{+}\triangleright(m\ot_{\mathcal{K}} t_0\ns{0}t_1\ns{0}\ns{0}\cdots t_n\ns{0}\ns{0})\\
&=(t_0\ns{1}t_1\ns{1}\cdots t_n\ns{1}\ot_R t_1\ns{2}\cdots t_n\ns{2}\\
&\ot_R \cdots \ot_Rt_{n-1}\ns{n}t_n\ns{n})\triangleleft {t_{n}\ns{n+1}}^{-}\ot_{R^{op}} \\
&{t_{n}\ns{n+1}}^{+}\triangleright(m\ot_{\mathcal{K}} t_0\ns{0}t_1\ns{0}\cdots t_n\ns{0})\\
&=(t_0\ns{1}t_1\ns{1}\cdots t_{n-1}\ns{1}\ot_R t_1\ns{2}\cdots t_{n-1}\ns{2}\ot_R\\
& \cdots \ot_Rt_{n-1}\ns{n})\triangleleft t_n\ns{1}{t_{n}\ns{2}}^{-}\ot_{R^{op}} \\
&{t_{n}\ns{2}}^{+}\triangleright(m\ot_{\mathcal{K}} t_0\ns{0}t_1\ns{0}\cdots t_n\ns{0})\\
&=t_0\ns{1}t_1\ns{1}\cdots t_{n-1}\ns{1}\ot_R t_1\ns{2}\cdots t_{n-1}\ns{2}\ot_R \cdots \ot_Rt_{n-1}\ns{n}\ot_{R^{op}} \\
&t_n\ns{1}\triangleright(m\ot_{\mathcal{K}} t_0\ns{0}t_1\ns{0}\cdots t_n\ns{0})\\
&=t_0\ns{1}t_1\ns{1}\cdots t_{n-1}\ns{1}\ot_R t_1\ns{2}\cdots t_{n-1}\ns{2}\ot_R \cdots \ot_Rt_{n-1}\ns{n}\ot_{R^{op}} \\
&m\ns{0}\ot_{\mathcal{K}} {t_n\ns{1}}\ns{+}(m\ns{-1}\triangleright (t_0\ns{0}t_1\ns{0}\cdots t_n\ns{0}{t_n\ns{1}}\ns{-})\\
&=t_0\ns{1}t_1\ns{1}\cdots t_{n-1}\ns{1}\ot_R t_1\ns{2}\cdots t_{n-1}\ns{2}\ot_R \cdots \ot_Rt_{n-1}\ns{n}\ot_{R^{op}} \\
&m\ns{0}\ot_{\mathcal{K}} t_n(m\ns{-n}\triangleright t_0\ns{0})(m\ns{-n+1}t_1\ns{0})\cdots (m\ns{-1}t_{n-1}\ns{0})\\
&=\omega_n(m\ns{0}\ot_{\mathcal{K}} t_n \ot_S m\ns{-n}t_0\ot_S\cdots \ot_S m\ns{-1}t_{n-1})\\
&=\omega_n\tau_n(m\ot_{\mathcal{K}} t_0\ot_S\cdots \ot_S t_n).
  \end{align*}
We use the coaction property of $T$ over $\mathcal{B}$ in the second equality, the right diagonal action of $\mathcal{B}$ over $\mathcal{B}^{\ot_R n}$ in the third equality, Lemma \cite[2.14 .(ii)]{bs2} in the fourth equality, the action \eqref{action} in the fifth equality, Lemma \ref{inverse}(vi) in the sixth equality, the $\mathcal{K}$-module algebra property of $T$ in the seventh equality and $\mathcal{K}$-equivariant property of the coaction of $\mathcal{B}$ over $T$ in the eighth equality.
\end{proof}
%%%%%%%%%%%%%%%%%%%%%%%%%%%%%%%%%%%%%%%%%%%%%%%%%%%%%%%%%%%%%%%%%%%%%%%%%%%%%%%%%%%%%%%%%%%%%%%%%%%%%%%%%%%%%%%%%%%%%%%%%%%%%%%%%%%%

\subsection{Examples}\label{ss-3-2}
As the first example,  one notes  that by  using  Proposition \ref{nice} and Theorem\ref{main},  for $\mathcal{K}= S\ot S^{op}$ and $M= S$,  we see that    Theorem \ref{general} generalizes Theorem 2.6 at \cite{bs2}.

We like to illustrate the SAYD module constructed in Theorem \ref{main} in  an explicit and nontrivial example.
Let $\Hc$ and $\Fc$  denote  two Hopf algebras such that $\Hc$ acts on  $\Fc$ and makes it a left module algebra. We let $A:=\Fc\al \Hc$ be the usual crossed product algebra, that is $\Fc\ot \Hc$ as a vector space  with multiplication $(f\al h)(g\al v)= f(h\ps{1}\rt g)\al h\ps{2}v$, and $1\al 1$ as its unit. One knows that $\Db: A\ra A\ot \Hc$ defined by $\Db(f\al h)=f\al h\ps{1}\ot h\ps{2}$ defines a comodule algebra over the Hopf algebra $\Hc$. One also easily verifies that $B:=A^\Hc=\Fc\ot \Cb\simeq \Fc$, and $A(B)^\Hc$ is  a Hopf Galois extension. We let $\Kc:= B\ot \Fc\ot B^{\rm o}$ be the usual left $\times_B$-Hopf algebra,  i.e., as an algebra it is tensor product of algebras and its coring structure is given  by,
\begin{align}
\begin{split}
&\mathfrak{s}(b)= b\ot 1\ot 1, \quad \mathfrak{t}(b)= 1\ot 1\ot b,\\
&\D(b_1\ot f\ot b_2)= b_1\ot f\ps{1}\ot 1\ot_B 1\ot f\ps{2}\ot b_2, \\
& \ve(b_1\ot f\ot b_2)= \ve(f)b_1b_2,
\end{split}
\end{align}
One observes that
\begin{equation}
\nu^{-1}( b_1\ot f\ot b_2\ot _{B^o} b_3\ot g\ot b_4)= b_1\ot f\ps{1}\ot 1\ot_B b_2b_3\ot S(f\ps{2})g\ot b_4.
\end{equation}
 So we see that
 \begin{equation}
 (b_1\ot f\ot b_2)^-\ot_{B} (b_1\ot f\ot b_2)^+= b_1\ot f\ps{1}\ot 1\ot_B b_2\ot S(f\ps{2})\ot 1.
 \end{equation}

Now let $N$ be a right-left SAYD module over the Hopf algebra $\Fc$ and let also  $x\ot x^{-1}\in B\ot B^o$ be a group-like element with $x\in Z(B)$, the center of $B$.
 We  define the following action and coaction on   $M:= ~^{x\ot x^{-1}}B\ot N=B\ot N$.
\begin{align}
&\Db(b\ot n)= bx\ot n\ns{-1}\ot x^{-1}\ot_B 1\ot n\ns{0}, \\
&(b\ot n)\cdot (b_1\ot f\ot b_2)= b_2bb_1\ot n\cdot f
 \end{align}

\begin{lemma}
Via the above action and coaction $~^{x\ot x^{-1}}B\ot N$ is a right-left SAYD module over $\Kc$.
\end{lemma}
\begin{proof}
It is obvious that the action and coaction are well-defined.  In the following we check the AYD condition,
\begin{align*}
&\Db((b\ot n)\cdot (b_1\ot f\ot b_2))= \Db(b_2bb_1\ot n\cdot f)=\\
&b_2bb_1 x\ot (nf)\ns{-1}\ot x^{-1}\ot_B 1\ot (n\cdot f)\ns{0}=\\
&b_2bb_1x\ot S(f\ps{3})n\ns{-1}f\ps{1}\ot x^{-1}\ot_B 1\ot  n\ns{0}\cdot f\ps{2}=\\
&(b_2\ot S(f\ps{3})\ot 1)(bx\ot n\ns{-1}\ot x^{-1})(b_1\ot f\ps{1}\ot 1)\ot_B (1\ot n\ns{0})\cdot (1\ot f\ps{2}\ot 1)=\\
&{k\ps{2}}^+n\ns{-1} {k\ps{1}}\ot_B  n\ns{0}\cdot {k\ps{2}}^-\,.
\end{align*}
The stability condition follows from the stability of $N$.
\end{proof}
One then defines a left action of $\Kc$ on $A$ by
\begin{equation}\label{K-acts-A}
(b_1\ot f\ot b_2)\rt (g\al h)= b_1 f\ps{1} g \;(h\ps{1}\rt (S(f\ps{2}) b_2))\al h\ps{2}.
\end{equation}
\begin{lemma}
Let $\Fc$ be commutative. Then   the above action makes $A$ a $\Kc$-module algebra.
\end{lemma}
\begin{proof}
Let us first prove that the equation \eqref{K-acts-A} defines an action.
Indeed for $k_1:= b_1\ot e\ot b_2$, and $k_2:= b_3\ot f\ot b_4$ by using the facts that $\Fc$ is commutative and $\Fc$ is $\Hc$-module algebra we have
\begin{align*}
&k_1\rt (k_2\rt (g\al h))= k_1\rt (b_3 f\ps{1} g(h\ps{1}\rt (S(f\ps{2}) b_4))\al h\ps{2})=\\
& (b_1b_3 e\ps{1}f\ps{1} g(h\ps{1}\rt (S(f\ps{2} b_4)  h\ps{2}\rt(S(e_2)b_2)\al h\ps{2})=k_1k_2\rt(g\al h).
\end{align*}
Obviously the action defined in \eqref{K-acts-A} is unital.
We use $h\rt 1_\Fc=\ve(h)1_\Fc$ to see that
\begin{align*}
&\mathfrak{s}(b)\rt (f\al h)= (b\ot 1\ot 1)\rt (f\al h)= bf\al h. \\
&\mathfrak{t}(b)\rt (f\al h)= (1\ot 1\ot b)\rt (g\al v)= g (h\ps{1}\rt b)\al h\ps{2}.
\end{align*}
This tells us that the multiplication of $\Kc$ is $B$-balanced.
Now we show that $k\rt (a_1a_2)=(k\ps{1}\rt a_1)(k\ps{2}\rt a_2)$. Let $k= (b_1\ot f\ot b_2)$,
 $a_1= g\al h$, and $a_2= l\al v$, we  have
\begin{equation*}
k\rt (a_1a_2)= k\rt (g h\ps{1}\rt l\al h\ps{2}v)= b_1 f\ps{1}g (h\ps{1}\rt l) (h\ps{2}v\ps{1}\rt (S(f\ps{2}) b_2))\al h\ps{3}v\ps{2}.
\end{equation*}
On the other hand by using the fact that $\Fc$ is $\Hc$-module algebra we see
\begin{align*}
&(k\ps{1}\rt a_1)(k\ps{2}\rt a_2)= ((b_1\ot f\ps{1}\ot 1)\rt (g\al h))((1\ot f\ps{2}\ot b_2)\rt (l\al v))=\\
& (b_1 f\ps{1} g \,(h\ps{1}\rt S(f\ps{2}))\al h\ps{2})(f\ps{3} l (v\ps{1}\rt (S(f\ps{4})b_2))\al v\ps{2})=\\
& b_1 f\ps{1} g \,(h\ps{1}\rt S(f\ps{2}))(h\ps{2}\rt (f\ps{3} l (v\ps{1}\rt (S(f\ps{4})b_2)))\al h\ps{3}v\ps{2}=\\
&b_1 f\ps{1} g\, (h\ps{1}\rt S(f\ps{2}))(h\ps{2}\rt f\ps{3}) (h\ps{3}\rt l) (h\ps{4}v\ps{1}\rt (S(f\ps{4})b_2))\al h\ps{5}v\ps{2}=\\
&b_1 f\ps{1}g\, (h\ps{1}\rt l) (h\ps{2}v\ps{1}\rt (S(f\ps{2}) b_2))\al h\ps{3}v\ps{2}.
\end{align*}
Finally we check the condition \eqref{modulealgebra-1}. For $k:= b_1\ot f\ot b_2$   we have,
\begin{align*}
&k\rt 1_A= (b_1\ot f\ot b_2)\rt (1\al 1)= b_1f\ps{1}(S(f\ps{2})b_2)\al 1=\\
&\ve(f)b_1b_2\al 1= (\ve(f)b_1b_2\ot 1\ot 1)\rt (1\al 1)=\Fs(\ve(k))\rt 1_A.
\end{align*}
 \end{proof}
 One easily sees that the coaction $\Db: A\ra A\ot \Hc$ is $\Kc$ equivariant.

\bigskip

 We now let the group of all unit elements of $B$, which is denoted by $B^{\times}$, define a category whose objects and morphisms are both elements of $B^\times$. More precisely, for $b_1,b_2\in B^\times$,  ${\rm Mor}(b_1,b_2)$ consists the unique element denoted by $\hat c$, where  $c=b_2b_1^{-1}$.  We let  $B^{\times}\times~ ^\Fc SAYD_\Fc$ be the product of categories. Here  $^\Fc SAYD_\Fc$  is the category  whose objects are  left-right SAYD modules over $\Fc$ and morphisms are $\Fc$-linear and $\Fc$-colinear maps.
 We  define the following functor
\begin{align}
\begin{split}
&\Phi: B^{\times}\times~ ^\Fc SAYD_\Fc\ra ~_\Kc SAYD(\Hc)^\Kc,\\
&\Phi((x,N))=~~ ^{x\ot {x^{-1}}}B\ot N,\\
&\Phi(\hat y,\phi) : ~ ^{x\ot {x^{-1}}}B\ot N_1\ra ^{xy\ot {(xy)^{-1}}}B\ot N_2 ,\\
&  \Phi(\hat y,\phi) (b\ot n)=by^{-1}\ot \phi(n).
\end{split}
 \end{align}
 Here    $ y\in  B^\times $ and $\phi\in~ ^\Fc\Hom_\Fc(N_1,N_2)$.

 \begin{proposition}\label{proposition-Phi-functor}
 The above assignment $\Phi$ defines a covariant functor.
 \end{proposition}
 \begin{proof}
 We need to show that $\Phi$ is a morphism map.  Indeed, first we see that $\Phi(\hat y,\phi)$ is a $\Kc$-colinear map. Using the facts that $\phi$ is a $\Fc$-colinear map, that $B$ is commutative,   and also the bialgebroid structure of $\Kc=B\ot \Fc\ot B$,  we see
 \begin{align*}
  &\Db(\Phi(\hat y,\phi)(b\ot n))= \Db( by^{-1},\phi(n))=\\
 &by^{-1}xy\ot \phi(n)\ns{-1}\ot (xy)^{-1}\ot_B 1\ot \phi(n)\ns{0}=\\
 &bx\ot \phi(n\ns{-1})\ot x^{-1}y^{-1}\ot_B 1\ot \phi(n\ns{0})=\\
 &bx\ot \phi(n\ns{-1})\ot x^{-1}\ot_B y^{-1}\ot \phi(n\ns{0})=\\
 &( \Id_\Kc \ot \Phi(\hat y,\phi))(\Db(b\ot n)).
 \end{align*}
 Now we prove that $\Phi(\hat y,\phi)$ is a $\Kc$-linear map.
 \begin{align*}
 &\Phi(\hat y,\phi)((b\ot n)\cdot(b_1\ot f\ot b_2))= \Phi(\hat y,\phi)(b_2bb_1\ot n\cdot f)=\\
 &(b_2bb_1y^{-1}\ot \phi(n\cdot f))=(b_2bb_1y^{-1}\ot \phi(n)\cdot f)= \\
 &(b_2bb_1y^{-1}\ot \phi(n\cdot f)) =(\Phi(\hat y,\phi))(b_2bb_1\ot n\cdot f).
 \end{align*}

 Finally one uses again the commutativity of $B$ to see that $\Phi((\hat{y_1},\phi_1)\circ (\hat{y_2},\phi_2))=\Phi(\hat{\overbrace{y_1y_2}}, \phi_1\phi_2)= \Phi(\hat{y_1}, \phi_1)\circ \Phi(\hat{y_2}, \phi_2)$ .
 \end{proof}

 As a result one composes the functors $\Phi$ and  $\digamma$ defined in Propositions \ref{diagmma} and \ref {proposition-Phi-functor} respectively to get the following functor
  \begin{equation}
\digamma \circ \Phi: B^{\times}\times~ ^\Fc SAYD_\Fc\ra ~_\Hc SAYD^\Hc.
 \end{equation}
One notes that in the simplest possible case of this example, i.e, when $\Fc:=B=\Cb$,  $N:=\Cb$, $x=1_\Cb$ the resulting SAYD module is $\Hc$ with the  standard action and coaction, i.e, the action and coaction are  defined by the  adjoint action and  the  comultiplication of $\Hc$ respectively.

%%%%%%%%%%%%%%%%%%%%%%%%%%%%%%%%%%%%%%%%%%%%%%%%%%%%%%%%%%%%%%%%%%%%%%%%%%%%%%%%%%%%%%%%%%%%%%%%%%%%%%%%%%%%%%%%%%%%%%%%%

%%%%%%%%%%%%%%%%%%%%%%%%%%%%%%%%%%%%%%%%%%%%%%%%%%%%%%%%%%%%%%%%%%%%%%%%%%%%%%%%%%%%%%


\begin{thebibliography}{9}


%%%%%%%%%%%%%%%%%%%

\bibitem[B]{b1}
G. B\"{o}hm, {Hopf Algebroids}, Handbook of Algebra Vol \textbf{6}, edited by M. Hazewinkel, Elsevier( 2009),  pp. 173--236.

\bibitem[BS]{bs2}
G. B\"{o}hm, and D. Stefan, {(co)cyclic (co)homology of bialgebroids: An approach via (co)monads}, Commun. Math. Phys\textbf{ 282} (2008),  pp. 239--286.
\bibitem[BW]{bwbook}
T. Brzezinski, and  R. Wisbauer, {Corings and Comodules},  London Mathematical Society Lecture Note Series \textbf{ 309}
Cambridge University Press, Cambridge (2003).

\bibitem[C-Book]{NCG} A. Connes, {\bf Noncommutative geometry},  Academic Press, 1994.

\bibitem[CM98]{ConMos:HopfCyc}  A. Connes,  and H. Moscovici, \emph{Hopf algebras,
  cyclic cohomology and the transverse index theorem},
  Comm. Math. Phys. \textbf{198} (1998), pp. 199--246.

\bibitem[CM01]{ConMos:DiffCyc} A. Connes,   and H. Moscovici, \emph{Differential
    cyclic cohomology and Hopf algebraic structures in transverse geometry},
    in: Essays on geometry and related topics. Vol 1-2,
    Monogr. Enseign. Math. \textbf{ 38}, Enseignement Math, Geneva (2001), pp. 217--255.


  \bibitem[CM00]{ConMos:Cyc-Symm} A.  Connes,  and Moscovici, H.: Cyclic cohomology and Hopf algebra symmetry, Lett. Phys. 52 (1), (2000), 1�28.

\bibitem[Cra]{Cra} M. Crainic, {\em Cyclic cohomology of \'etale groupoids: the
    general case}, $K$-Theory \textbf{17} (1999), 319--362.

\bibitem[Gor]{Gor:SecChar}
Gorokhovsky, A., {Secondary characteristic classes and cyclic
cohomology of Hopf algebras}, {\it Topology} {\bf 41} (2002), 993-1016.


\bibitem[HKRS1]{HaKhRaSo1} P.M. Hajac, M. Khalkhali, B. Rangipour and
  Y. Sommerh\"auser, \emph{Stable anti-Yetter-Drinfeld modules},
  C. R. Math. Acad. Sci. Paris  \textbf{338}  (2004), pp. 587--590.

\bibitem[HKRS2]{HaKhRaSo2} P.M. Hajac, M. Khalkhali, B. Rangipour and
  Y. Sommerh\"auser, \emph{Hopf-cyclic homology and cohomology with
  coefficients},
  C. R. Math. Acad. Sci. Paris  \textbf{338}  (2004), pp. 667--672.


\bibitem[JS]{JaSt:CycHom} P. Jara and D. Stefan, \emph{Hopf-cyclic
    homology and relative cyclic homology of Hopf-Galois extensions},
    Proc. London Math. Soc. (3) \textbf{93} (2006), pp. 138--174.

\bibitem[Kay05]{Kay} A. Kaygun, \emph{Bialgebra cyclic homology with
  coefficients}, $K$-Theory \textbf{34} (2005), pp. 151--194.

\bibitem[Kay06]{Kay:UniHCyc} A. Kaygun, \emph{The universal Hopf cyclic
  theory},  Journal of Noncommutative Geometry Vol.\textbf{2} (2008), No.2, pp. 333�-351

\bibitem[KP]{KhaPou} M. Khalkhali, and  A. Pourkia, \emph{ Hopf cyclic cohomology in braided monoidal categories.}  Homology, Homotopy Appl.  {\bf 12}  (2010),  no. 1, 111--155.



\bibitem[KR05]{KhaRan:ANoteCycDual} M. Khalkhali and B. Rangipour, \emph{A note on cyclic
    duality and Hopf algebras},  Comm. Algebra \textbf{33} (2005), pp. 763--773.


\bibitem[KR04]{KhaRan:Para_Hopf} M. Khalkhali and B. Rangipour, \emph{Para-Hopf
    algebroids and their cyclic cohomology}, Lett. Math. Phys. \textbf{70}
    (2004), pp. 259--272.
\bibitem[KR03]{KhaRan:InvInvCycHom} M. Khalkhali and B. Rangipour,
  \emph{Invariant cyclic homology}, $K$-Theory \textbf{28} (2003), pp. 183--205.
\bibitem[KR02]{KhaRan:ANewCyc}
M. Khalkhali and B. Rangipour, A new cyclic module for Hopf algebras,  K-Theory  27  (2002), pp. 111--131.





\bibitem[Lo]{Loday} J.-L. Loday, {\em Cyclic Homology}. Springer-Verlag 1992.



\bibitem[Ra]{Rang:CycCor} B. Rangipour, {\em Cyclic cohomology of corings},
Journal of K-theory, Cambridge University Press (2009), \textbf{4}, pp. 193--207.

\bibitem[Sch98]{Scha:bia_nc} P. Schauenburg, {\em Bialgebras over
  noncommutative rings and a structure theorem for Hopf bimodules},
  Appl. Categorical Str. \textbf{6} (1998), pp. 193--222.

\bibitem[Sch]{sch} P. Schauenburg, {\em Duals and doubles of
  quantum groupoids ($\times_R$-Hopf algebras)}, in: N. Andruskiewitsch,
  W.R. Ferrer-Santos and H.-J. Schneider (eds.) AMS
  Contemp. Math. \textbf{267}, AMS Providence (2000) pp. 273--293.

\end{thebibliography}
  \end{document}